\documentclass[preprint,12pt]{elsarticle}

\usepackage{amssymb}
\usepackage{mathtools}
\usepackage{amsmath}
\usepackage{cases}
\usepackage{float}
\usepackage{comment}
\usepackage[unicode=true]{hyperref}
\journal{Communications in Nonlinear Science and Numerical Simulations }

\begin{document}

\begin{frontmatter}


\title{Mathematical models for CAR T immunotherapy and CD19 dynamics in leukemia: a comparative analysis} 

\author[af1,af2]{Salvador Chulián\corref{cor1}}
\ead{salvador.chulian@uca.es}
\author[af1,af2,af3]{Ana Niño-López}
\author[af1,af2]{Rocío Picón-González} 
\author[af1,af2]{María Rosa} 
\affiliation[af1]{organization={Department of Mathematics, University of Cádiz},city={Puerto Real},country={Spain}}
\affiliation[af2]{organization={Biomedical Research and Innovation Institute of Cádiz (INiBICA), Puerta del Mar University Hospital},city={Cádiz},country={Spain}}
\affiliation[af3]{organization={Pediatric Health Research Institute Niño Jesús University Children's Hospital (IPIS-NJ)}, city={Madrid}, country={Spain}}
\begin{abstract}

Chimeric Antigen Receptor (CAR) T cell therapy has emerged as a successful treatment for relapsed or refractory hematological malignancies, particularly for B cell Acute Lymphoblastic Leukemia (B ALL), where CD19 targeted therapies have achieved high initial remission rates. However, relapse after treatment remains a major clinical challenge, frequently associated with antigen escape mechanisms and the emergence of CD19$^-$ leukemic cells. Understanding the interaction between CAR T cells and antigen expression dynamics is, therefore, essential for improving therapeutic efficacy and long-term patient outcomes.

In this work, we develop and comparatively analyze novel mathematical models describing CAR T immunotherapy and CD19 dynamics in leukemia. Our proposed framework combines compartmental ordinary differential equation (ODE) formulations as well as partial differential equation (PDE) systems. Both methods are able to capture the evolution of leukemic populations under immune pressure and, in particular, the models explicitly distinguish between CD19$^+$ and CD19$^-$ leukemic cells and incorporate bidirectional phenotypic transitions regulated by CAR T activity. 

The developed models provide biologically interpretable and computationally efficient tools for studying treatment response, resistance, and relapse mechanisms in CAR T cell therapy. We compare the ability of the different modeling approaches to reproduce CD19 antigen modulation and CAR T efficacy dynamics. We also propose sensitivity analyses to study the parameters' influence on the models' dynamics. Our work contributes to the mathematical understanding of antigen-driven resistance and offers a basis for future optimization and personalization of CAR T therapeutic strategies in leukemia.
\end{abstract}

\begin{graphicalabstract}

\includegraphics[width=\textwidth]{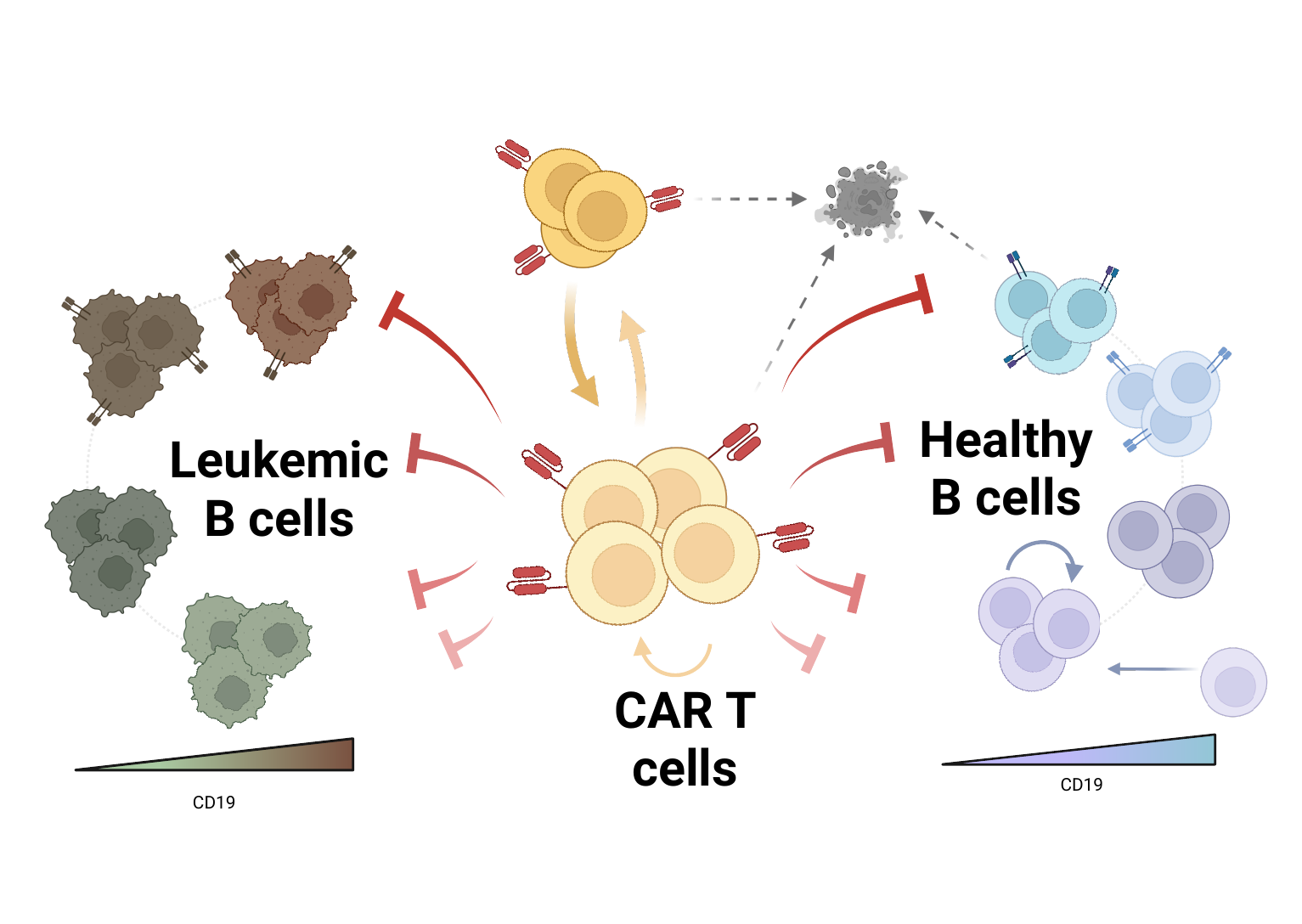}
\end{graphicalabstract}

\begin{highlights}
\item Novel ODE and PDE frameworks model CAR T cell immunotherapy and CD19 dynamics. 
\item Models capture bidirectional, immune-regulated phenotypic transitions in leukemia.
\item PDE system continuously structures leukemic populations by CD19 expression levels. 
\item Sensitivity analyses identify key parameters governing antigen escape and relapse.
\end{highlights}

\end{frontmatter}



\section{Introduction}
\label{sec1}

In recent years, Chimeric Antigen Receptor (CAR) T cell therapy has motivated numerous medical and mathematical studies aimed at evaluating its efficacy and optimizing its application in solid tumors, leukemias, and lymphomas \cite{globerson2021car, hay2017chimeric, swanson2022mathematical, bock2025outcome}. This immunotherapeutic approach has revolutionized the treatment of relapsed or refractory cancers, impacting both tumor evolution and current therapeutic strategies. In particular, in B cell Acute Lymphoblastic Leukemia (B ALL), CAR T therapy targeting the CD19 surface antigen has achieved response rates exceeding $80\%$ in patients who relapse after conventional treatments. However, despite these high initial response rates in both pediatric and adult populations, post-treatment relapse remains one of the major clinical challenges. Among the most relevant resistance mechanisms is antigen loss, especially the emergence of CD19$^{-}$ leukemic cells capable of evading CAR T-mediated immune activity.\\

Several clinical studies have highlighted important limitations associated with both CAR T cell persistence and tumor escape mechanisms driven by the loss or downregulation of CD19 expression. CD19$^+$ (CD19$^+$) and CD19$^-$ (CD19$^-$) relapses following different treatments have been analyzed in patient cohorts \cite{dourthe2021determinants}, while additional studies have investigated how pretherapy factors, including baseline CD19 expression, the presence of CD19$^-$ subpopulations, and prior anti CD19 therapies, influence clinical response \cite{zhang2020efficacy, pillai2019car}. Other approaches have explored strategies to improve therapeutic efficacy through CAR T reinfusion protocols \cite{gauthier2021factors}, as well as multitarget therapies incorporating additional antigens such as CD20 and CD22 to reduce immune escape phenomena \cite{fousek2021car}.\\

In parallel, mathematical modeling has played a crucial role in understanding these phenomena, ranging from early pharmacokinetic and cellular interaction models to more sophisticated formulations focused on immune escape mechanisms \cite{serrano2024understanding, sabir2025mathematical}. Several models have described CAR T expansion and B cell aplasia dynamics \cite{leon2021car, martinezrubio2021}, while other approaches have addressed T cell Acute Lymphoblastic Leukemia specific challenges such as fratricide \cite{perez2021car}. However, most existing models focus exclusively on CD19$^{+}$ relapses or treat CD19$^-$ states as absorbing states, without considering the possibility of antigen re expression or reversible phenotypic transitions.\\

More recently, several studies have emphasized the importance of mathematical and computational methodologies for refining CAR T therapeutic strategies and predicting treatment efficacy \cite{nukala2021systematic, putignano2025mathematical}. Mathematical models have been developed to describe CAR T therapy dynamics in B ALL and to parameterize cellular transition processes for patient stratification and minimal residual disease detection after chemotherapy \cite{haries2025effectiveness, gravenmier2026cell}. As well, artificial intelligence based approaches have shown considerable potential for improving minimal residual disease detection \cite{seheult2026artificial}, while in vitro experimental studies have enabled the estimation of key parameters required to calibrate dynamic interaction models between CAR T cells and leukemic cells \cite{shah2026data}.\\

Motivated by these limitations and recent advances, we introduce two novel compartmental frameworks that distinguish between CD19$^+$ and CD19$^-$ leukemic populations while allowing immune regulated bidirectional transitions between both states. Transition rates are modeled using Michaelis Menten kinetics \cite{meral2025mathematical} and are explicitly driven by CAR T cell activity, thereby capturing in a unified manner both immune escape and tumor resensitization dynamics, this is, how to make tumor cells responsive to therapies that previously failed. The resulting model remains within an ordinary differential equation (ODE) setting, depending only on time, which preserves biological interpretability, enables efficient computation, and facilitates calibration for future patient level clinical data. In addition, we extend this framework by proposing a partial differential equation (PDE) model that incorporates a continuous description of CD19 expression, allowing a more detailed representation of antigen heterogeneity over time. By incorporating reversible antigen modulation, this work goes beyond the assumptions of current approaches and provides a more mechanistic description of relapse. Ultimately, it offers a flexible and interpretable tool for dissecting resistance pathways and optimizing CAR T based therapeutic strategies in B ALL.

\section{Mathematical models}
In this section, we present two mathematical models that represent the dynamics of immunotherapy with leukemic and healthy cell populations. 

Generally, we have several subpopulations: healthy B cells, leukemic cells and CAR T cells. Depending on the antigen levels in both leukemic and healthy B cells, CAR T cells can activate themselves, react and kill those cells expressing such antigen. In this work, we are presenting two different mathematical models: In a first model we consider an ordinary differential system of equations (ODE model) depending on time, that compartmentalize leukemic cells with either a positive or a negative expression of the antigen; in a second model, we derive a partial differential equation system (PDE model) where the cell populations depend continuously on the antigen level as well as time. 

Specifically, we consider that leukemic cells can have either positive or negative CD19 antigen expression. We have thus proposed that the consideration of such levels can be a binary variable (either CD19$^+$ or CD19$^-$ for the ODE model), or that such antigen expression can be represented as a continuous variable (in the PDE model).

\subsection{ODE model}
\label{ODE model}
\begin{figure}[H]
    \centering
    \includegraphics[width=\linewidth]{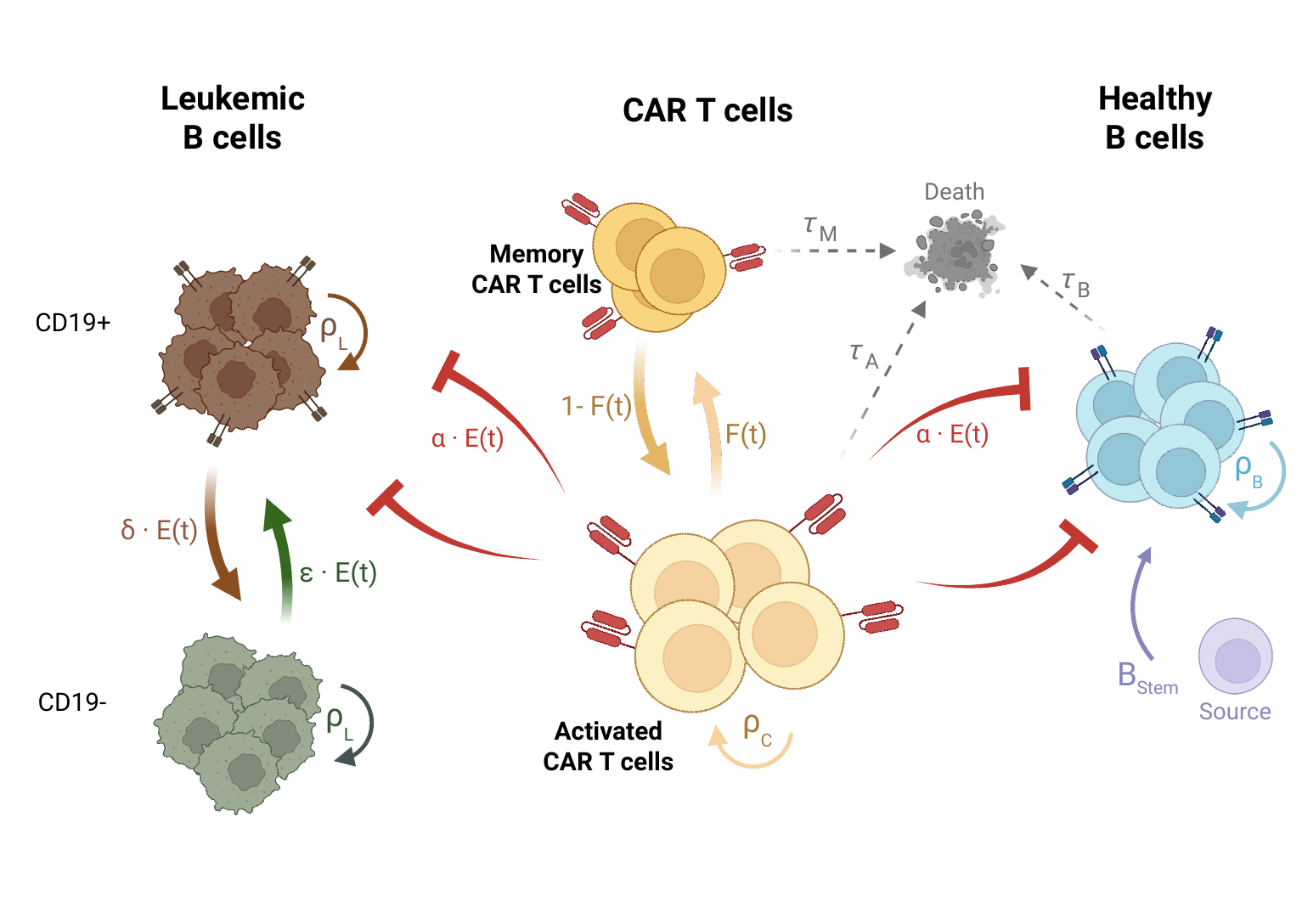}
    \caption{\textbf{Diagram of the ODE Model from Eqs. \eqref{EQ_ODE_model}}. In the ODE model, Leukemic B cells grow with rate $\rho_L$ and are divided in leukemic antigen positive cells ($L_P$) and antigen negative cells ($L_N$), depending on the therapy effect function $E(t)$ and parameters $\delta$ and $\varepsilon$. Healthy B cells ($B$) arise from a source $B_\text{stem}$, with growth $\rho_B$ and death rate $1/\tau_{B}$. CAR T cells are divided in memory and activated CAR T cells and their dynamics depend on the threshold activation function $F(t)$, dying with respective death rates $1/\tau_{M}$ and $1/\tau_{A}$. Activated CAR T cells kill antigen positive cells, both leukemic and healthy, depending on the effect of therapy $E(t)$ with rate $\alpha$.}
    \label{fig:ODEmodel}
\end{figure}

The ODE model, expressed in Eqs. \eqref{EQ_ODE_model} and where leukemic cells have a binary expression of the antigen (either CD19$^+$ or CD19$^-$), has been considered in Figure \ref{fig:ODEmodel}.

The system of Eqs. \eqref{EQ_ODE_model} for the ODE model  describes the interactions between healthy B cells $(B=B(t))$, leukemic B cells with and without CD19 antigen $(L_P=L_P(t), L_N=L_N(t))$, and CAR T cells in both their activated $(C_A=C_A(t))$ and memory states $(C_M=C_M(t))$. Each equation contains terms representing biological processes relevant to B cell leukemia and CAR T therapy:

\begin{equation}
\label{EQ_ODE_model}
\hspace{-1cm}
\left\{\begin{array}{cll}

\dfrac{dB(t)}{dt} =& B_{\text{Stem}} + \rho_B B(t)\left(1-\dfrac{B(t)+L_P(t)+L_N(t)}{B_{\text{max}}+L_P(t)+L_N(t)}\right) - \dfrac{B(t)}{\tau_B} 
- \alpha E(t) B(t), \vspace{7pt}\\

\dfrac{dL_P(t)}{dt}=&\rho_L\,L_P(t)\,\left(1-\dfrac{L_P(t)+L_N(t)}{L_{\max}}\right)-\left(\alpha+\delta\right)E(t)L_P(t)+\epsilon\left(1-E(t)\right) L_N(t),\vspace{7pt}\\

\dfrac{dL_N(t)}{dt} =& \rho_L L_N(t)\left(1 - \dfrac{L_P(t)+L_N(t)}{L_{\max}}\right) 
+  \delta E(t)L_P(t) 
- \epsilon\left(1 - E(t)\right) L_N(t), \vspace{7pt}\\

\dfrac{dC_A(t)}{dt} =& F(t)\left(\rho_{C_A} C_A(t) + \gamma_{MA} C_M(t)\right) 
- \dfrac{C_A(t)}{\tau_{A}} - \gamma_{AM}(1 - F(t))C_A(t), \vspace{7pt}\\

\dfrac{dC_M(t)}{dt} =&\gamma_{AM} (1 - F(t))C_A(t) 
- \dfrac{C_M(t)}{\tau_{M}} - \gamma_{MA}F(t)C_M(t),\vspace{7pt}\\
E(t)=&C_A(t)\dfrac{B(t)+L_P(t)}{h+B(t)+L_P(t)},\vspace{7pt}\\
F(t)=&\dfrac{B(t)+L_P(t)}{k+B(t)+L_P(t)},
\end{array}\right.
\end{equation}

\noindent with initial conditions:

\begin{equation}
L_N(0) = {L_N}_0,\quad 
L_P(0) = {L_P}_0,\quad 
B(0) = B_0,\quad 
C_A(0) = {C_A}_0,\quad 
C_M(0) = {C_M}_0.
\end{equation}

\subsubsection*{Healthy B cells}
The equation for $B(t)$ includes four processes. The source term $B_{Stem}$ represents continuous input from bone marrow stem cells. The proliferation term $\rho_B B(t)$ accounts for natural expansion of mature B cells. We define a pseudologistic term that takes into consideration the total amount of leukemic cells as well as healthy B cells, so that leukemic cells are able to invade healthy B cells and thus stop their growth. Losses occur through natural death, modeled by $-B(t)/\tau_B$, where $\tau_B$ is the average lifespan. Finally, CAR T mediated killing is modeled by $-\alpha E(t)B(t)$, where $\alpha$ is the killing rate and $E(t)$ includes the activated CAR T cell action modulated by a term with Michaelis--Menten form that includes the action of activated CAR T cells and the amount of antigen present, this is, the total amount of $B$ and $L_P$ cells. 
\subsubsection*{Leukemic antigen-positive cells}
The equation for $L_P(t)$ combines intrinsic growth, CAR T pressure, and phenotypic switching. The growth term $\rho_L L_P(t)(1 - (L_P(t)+L_N(t))/L_{\max})$ models pseudologistic proliferation with carrying capacity $L_{\max}$, where the presence of a high amount of leukemic cells can influence each other. Activated CAR T interaction induces two effects: direct killing at rate $\alpha$, and conversion into antigen-negative cells at rate $\delta$, both scaled by the saturating function $E(t)$. In contrast, when activated CAR T stimulation is low, antigen-negative cells can revert into the antigen-positive pool, modeled by $\varepsilon(1-E(t)) L_N(t)$, where $\varepsilon$ is the reversion rate.

\subsubsection*{Leukemic antigen-negative cells}
The antigen-negative population $L_N(t)$ follows a similar structure. Cells proliferate pseudologistically at rate $\rho_L$. Under activated CAR T pressure, antigen-positive cells are converted into antigen-negative cells, expressed by $\delta E(t) L_P(t)$. When CAR T presence is low, antigen-negative cells are lost through the reverse transition to antigen-positive cells, considering the term $\varepsilon(1- E(t)) L_N(t)$.

\subsubsection*{Activated CAR T cells}
Activated CAR T $C_A(t)$ dynamics depend on antigen stimulation, encoded by $F(t) = \dfrac{B(t)+L_P(t)}{k+B(t)+L_P(t)}$, where $k$ is an activation threshold.
When stimulation is high ($k$ is low), $F(t)\approx1$, and activated CAR T expand at rate $\rho_{C}$ and memory cells are transformed into activated at rate $\gamma_{MA}$. When stimulation is low ($k$ is high), $F(t)\approx0$, and then activated CAR T cells do not grow that much and switch their phenotype to Memory CAR T cells. Activated CAR T are also removed by two mechanisms: natural death at rate $1/\tau_{A}$ and transition to memory phenotype at rate $\gamma_{AM}$ when antigen is scarce, this is, considering the term $(1-F(t))$.

\subsubsection*{Memory CAR T cells}
The memory pool $C_M(t)$ receives input from $C_A$ cells through the term $(1-F(t))\gamma_{AM}C_A(t)$, whenever antigen is absent or low. Memory cells die at rate $1/\tau_{M}$. Upon renewed antigen stimulation, they are converted back into activated CAR T at rate $\gamma_{MA}$, stimulated by the function $F(t)$.

\subsubsection*{CAR T Effect function}
The function $E(t)$ regulates therapy effect by coupling their dynamics to the amount of Activated CAR T cells. As such, the higher the amount of antigen levels, this is, the amount of $B$ and $L_P$ cells, the greater effect of the treatment. This is also regulated by the threshold $h$. The whole efficacy of the treatment would be considered as $\alpha E(t)$, where $\alpha$ is the killing rate.

\subsubsection*{Phenotypic switch between Activated and Memory CAR T cells function}
The function $F(t)$ regulates CAR T activation and memory transitions with an activation threshold $k$ and by coupling their dynamics to the antigen load. It increases with the abundance of healthy and leukemic antigen-positive B cells, approaching saturation when antigen levels are large. \\

Overall, the system captures key processes: continuous generation and turnover of healthy B cells, logistic growth and antigen-loss of leukemic cells, and dynamic activation, proliferation, and persistence of CAR T cells.

\subsection{PDE model}
\label{PDE model}

We will now describe with the system from Eqs. \eqref{EQ_PDE_model}, that describes interactions between healthy B cells, leukemic B cells, and CAR T cells structured by antigen expression $x$. In this PDE model, each equation includes terms representing key biological processes relevant to healthy B cell, leukemia B cells and CAR T therapy, as in Section \ref{ODE model}. This structured PDE system captures key mechanisms: continuous generation and turnover of healthy B cells ($B=B(t,x))$, antigen-dependent killing of both healthy and leukemic cells ($L=L(t,x))$, logistic proliferation for healthy B cells, exponential growth for leukemic cells, and antigen-stimulated proliferation and decay of CAR T cells (activated and memory CAR T cells, $C_A=C_A(t),C_M=C_M(t)$, respectively). The following system thus provides a framework to analyze how antigen heterogeneity influences therapeutic outcomes:

\begin{equation}
\hspace{-1cm}
\label{EQ_PDE_model}
\left\{\begin{array}{cl}
\dfrac{\partial B(t,x)}{\partial t}=& B_{\text{Stem}}S(x) + \rho_B B(t,x)\left(1-\frac{B(t,x)+\int_{x_{\min } }^{x_{\max }} L(t,x) dx}{K(x)+\int_{x_{\min } }^{x_{\max }} L(t,x) dx}\right) - \dfrac{B(t,x)}{\tau_B} 
- \alpha E(t,x) B(t,x), \vspace{7pt}\\

\dfrac{\partial L(t,x)}{\partial t}=&\rho_L\,L(t,x)\left(1-\frac{\int_{x_{\min } }^{x_{\max }} L(t,x) dx}{L_{\max}}\right)-\alpha E(t,x) L(t,x),\vspace{7pt}\\

\dfrac{dC_A(t)}{dt} =&F(t)\left(\rho_{C_A} C_A(t) + \gamma_{MA} C_M(t)\right) 
- \dfrac{C_A(t)}{\tau_{A}} - \gamma_{AM}(1 - F(t))C_A(t), \vspace{7pt}\\

\dfrac{dC_M(t)}{dt} =&\gamma_{AM} (1 - F(t))C_A(t) 
- \dfrac{C_M(t)}{\tau_{M}} - \gamma_{MA}F(t)C_M(t),\vspace{7pt}\\

S(x)=&\dfrac{s}{1+\exp(m(x-x_0))},\vspace{7pt}\\

E(t,x)=&C_A(t)\dfrac{x}{x+h},\vspace{7pt}\\

F(t)=&\dfrac{\int_{x_{\min } }^{x_{\max }} B(t)+L(t) dx}{k+\int_{x_{\min } }^{x_{\max }} B(t)+L(t)dx},\vspace{7pt}\\

\end{array}\right.
\end{equation}

\noindent with:

\begin{equation}
\begin{array}{cc}

L(0,x) = L_0\sum_{i=1}^n\mathcal{G}(x,\mu_i,\sigma_i),\quad 
K(x) = B_\text{max}\sum_{j=1}^m\mathcal{G}(x,\mu_j,\sigma_j),\quad \vspace{5pt}\\
B(0,x) = B_0\mathcal{G}(x,0,\sigma_B),\quad
C_A(0) = {C_A}_0,\quad 
C_M(0) = {C_M}_0.
\end{array}
\end{equation}

\begin{figure}
    \centering
    \includegraphics[width=\linewidth]{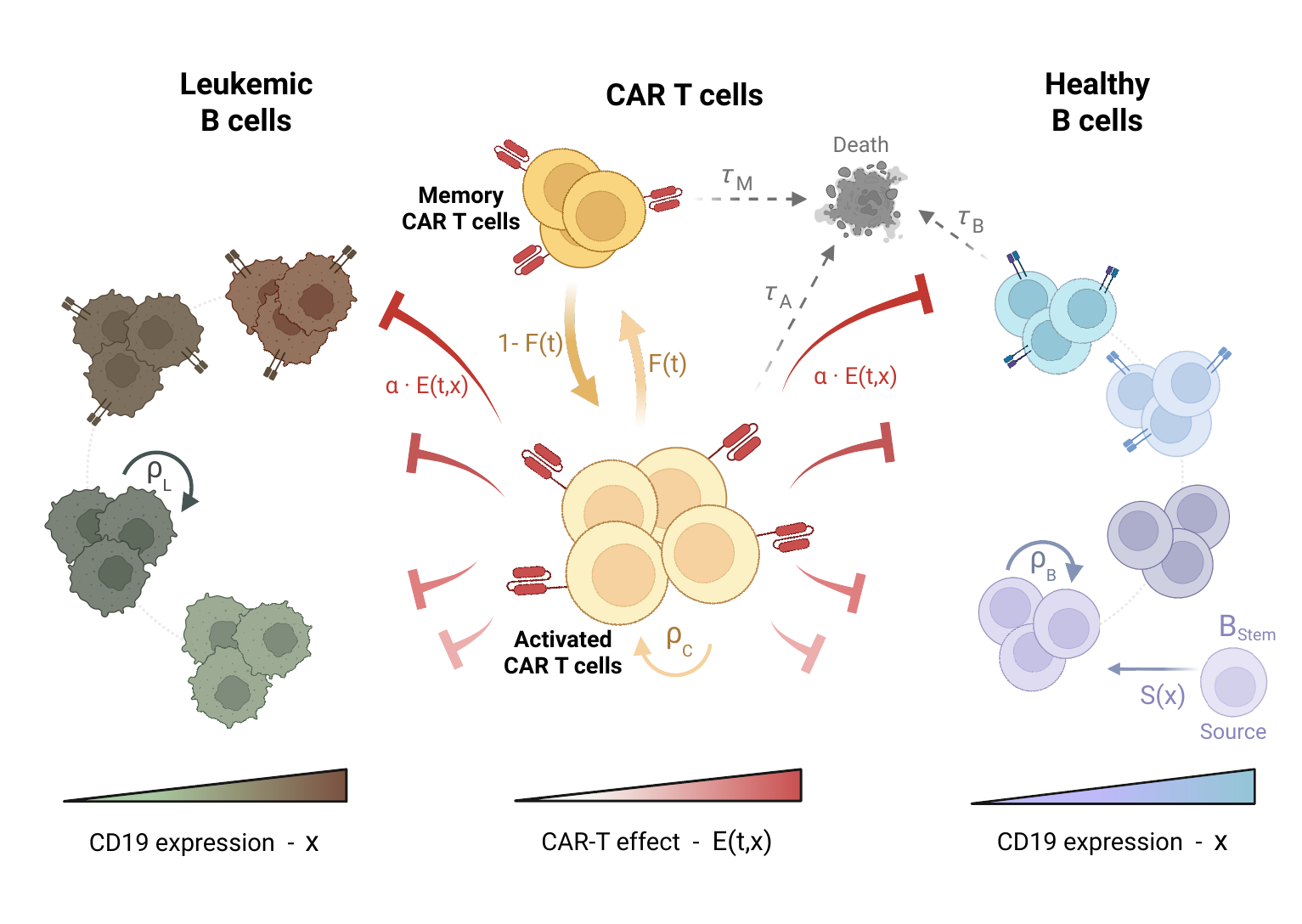}
    \caption{\textbf{Diagram of the PDE Model from Eqs. \eqref{EQ_PDE_model}}. In the PDE model, leukemic B cells grow with rate $\rho_L$, but also depend on the continuous expression of antigen $x$, here CD19. This antigen is expressed also in healthy B cells ($B$), which arise from a source $B_\text{Stem}$ (modulated by function $S(x)$), with growth $\rho_B$ and death rate $1/\tau_{B}$. CAR T cells, as in the previous ODE model, are divided in memory and activated CAR T cells depending on the threshold activation function $F(t)$, and die with respective death rates $1/\tau_{M}$ and $1/\tau_{A}$. The main difference for CAR T cells in comparison to the previous model, is that activated CAR T cells kill antigen positive cells, both leukemic and healthy, depending on the effect of therapy $E(t,x)$, which now depends on the antigen levels $x$,  with rate $\alpha$.}
    \label{fig:PDEmodel}
\end{figure}

We now describe each term in the following subsections:

\subsubsection*{Healthy B cells}

The equation for $B(t,x)$ includes four terms. The stem cell source term, $B_\text{Stem}S(x)$, models continuous input of a number of $B_\text{Stem}$, B cells from the bone marrow, biased toward cells with lower antigen expression with the function $S(x)$. The term $\rho_B B(t,x)\left(1-\frac{B(t,x)+\int_{x_{\min } }^{x_{\max }} L(t,x) dx}{K(x)+\int_{x_{\min } }^{x_{\max }} L(t,x) dx}\right)$ shows a proliferation rate $\rho_B$ and a pseudologistic growth term with antigen-dependent carrying capacity $K(x)$, which depends on antigen expression and is expressed as a sum of Gaussian distributions $\mathcal{G}$ centered at $\mu_j$ and standard deviation $\sigma_j$. Besides, the total amount of leukemic cells, expressed by the integral $\int_{x_{\min } }^{x_{\max }} L(t,x) dx$, can decelerate growth of healthy B cells due to invasion. Natural decay is expressed as $-B(t,x)/\tau_B$, representing average B cell lifespan. Finally, CAR T-mediated killing, $-\alpha E(t,x) B(t,x)$, depends on killing rate $\alpha$ and therapy effect $E(t,x)$, which now depends on antigen expression $x$.

\subsubsection*{Leukemic B cells}

Leukemic cells $L(t,x)$ grow intrinsically at rate $\rho_L$, with a logistic growth and assuming a general carrying capacity $L_\text{max}$ and the total amount of leukemic cells $\int_{x_{\min } }^{x_{\max }} L(t,x) dx$. They are eliminated by CAR T cells similarly as healthy B cells, with a killing rate $\alpha$ and therapy effect $E(t,x)$, reflecting antigen-dependent cytotoxicity. Unlike healthy B cells, leukemic cells in this model do not include stem cell input or natural cell death, this last reflecting aggressive growth behavior.

\subsubsection*{Sigmoid function for the influence of antigen into source term}
The function $S(x)$ models a sigmoid function with height $s$, that decays at a rate $m$ with center at $x_0$. This function helps on the biased action of the continuous input of $B_\text{Stem}$ cells for the healthy B cells, such that lower-antigen expressed cells have a higher demand of new cells coming from the bone marrow, in contrast to high-antigen expressed B cells.

\subsubsection*{CAR T cells, therapy effect and phenotypic switch function}

CAR T dynamics are driven similarly as in the previous ODE model. Nevertheless, CAR T Effect function $E(t,x)$ represent their corresponding therapy effect from activated CAR T cells, but now influenced by antigen expression: the higher the antigen expression $x$, the greater the therapy effect. Also, the phenotypic switch function $F(t)$ changes by including the total amount of $B(t,x)$ and $L(t,x)$, as an integration over the antigen dependent variable $x$ from a minimum and maximum point $x_\text{min}$ and $x_\text{max}$.

\subsubsection*{Initial conditions and function $K(x)$}

Populations $L(0,x)$ and $B(0,x)$ are initialized as sum of Gaussian distributions $\mathcal{G}(x,\mu_i,\sigma_i)$ over antigen $x$, representing natural variability in antigen expression for $B$ and $L$ cells. Carrying capacity function $K(x)$ is modeled also a sum of Gaussian functions, with a $B_\text{max}$ value.  The proportion of such cells in the sum of Gaussians depends on the different levels of antigen expression for healthy B cells \cite{CHULIAN2021110685}, where low, intermediate, and high-antigen cells have different proportions and are associated to different antigen expression levels $\mu_j$ with standard deviation $\sigma_j$. This is analogous for the different sets of leukemic clones, which can be simulated with several sums of Gaussian distributions centered at different antigen expression levels $\mu_i$ with standard deviation $\sigma_i$.

\section{Models' parameters and numerical implementation}

We present a general Table \ref{table param} of all parameters used and their corresponding values and meanings. Parameters source was included and typed as ``E'' if estimated range or ``A'' if it is an assumption.

\begin{table}[ht!]
\footnotesize{
\centering
\begin{tabular}{lllll}
\hline
\textbf{Param.} & \textbf{Meaning} & \textbf{Value} & \textbf{Units} & \textbf{Source} \\
\hline

\multicolumn{5}{l}{\quad\textit{\textbf{B cells}}} \\
\hline

$B_{\text{Stem}}$ &  Flux from Stem cells & $10^8$ & cell day$^{-1}$ & E \cite{martinezrubio2021,nino2023mathematical}\\ 
$\rho_{B}$ &  B cell proliferation rate &$\log(2)/8$ & day$^{-1}$ & \cite{leon2021car} \\ 
$B_{\text{max}}$ &  B cell carrying capacity (ODE) & $5\times 10^9$ & cell & E \cite{nino2023mathematical} \\ 
$\tau_{B}$ &  B cell lifetime & $60$ & day & \cite{leon2021car}\\
$s$ & Signal for influx intensity & $[10^{-3},10^{-1}]$ & -& E\\
$x_0$ & Antigen influx threshold &0.5 & -&A\\
$m$ & Decay slope for influx  &1&- &A\\
\hline

\multicolumn{5}{l}{\quad\textbf{\textit{L cells}}} \\
\hline

$\rho_{L}$ &  Leukemic B cell proliferation rate & $1/30$ & day$^{-1}$ & \cite{leon2021car} \\ 
$L_{\text{max}}$ &  Leukemic B cell carrying capacity & $10^{12}$ & cell & \cite{nino2023mathematical} \\ 
\hline

\multicolumn{5}{l}{\quad\textbf{\textit{CAR T cells}}} \\
\hline
$k$  &  CAR T activation threshold & [$10^{8},10^{12}$] &cell &\cite{martinezrubio2021}\\
$\rho_{CA}$ & Activated CAR T proliferation & $0.9$ & day$^{-1}$& \cite{mueller2018clinical,stein2019tisagenlecleucel} \\
$\tau_{A}$ & Activated CAR T lifetime & $6.5$ & day& \cite{mueller2018clinical,stein2019tisagenlecleucel} \\
$\gamma_{AM}$ & Activated to memory transition & $0.001$ & day$^{-1}$ & \cite{mueller2018clinical,stein2019tisagenlecleucel}\\
$\gamma_{MA}$ & Memory to activated transition & $0.33$ & day$^{-1}$ & \cite{martinezrubio2021}\\
$\tau_{M}$ & Memory CAR T lifetime & $300$ & day& \cite{mueller2018clinical,stein2019tisagenlecleucel} \\
$\alpha$ & Active CAR T killing capacity (PDE) & $3\times 10^{-10}$&cell$^{-1}$day$^{-1}$ & \cite{martinezrubio2021}\\
$h_O$  & Treatment efficacy activation threshold (ODE)  & $[10^8,10^{12}]$&cell & E\\
$h_P$  & Treatment efficacy activation threshold (PDE)  & [0,1] &- & E\\

\hline
\multicolumn{5}{l}{\quad \textbf{\textit{Initial Conditions and general parameters}}} \\
\hline
${C_A}_0$ & Initial active CAR T cells & $10^7$ & cell & E \cite{martinezrubio2021}\\
${C_M}_0$ & Initial memory CAR T cells & 0 & cell & \cite{martinezrubio2021}\\

${B}_0$ & Initial healthy B cells (ODE)&  $10^7$& cell & \cite{nino2023mathematical}\\
${L_P}_0$ & Initial leukemic antigen positive cells (ODE)&$10^{7}$ & cell & A\\
${L_N}_0$ & Initial leukemic antigen positive cells (ODE)& 0 or low value& cell & A\\
$x_{\min}$ & Minimal antigen expression $x$ &0 &- &A\\
$x_{\max}$ & Maximal antigen expression $x$ &1 &- &A\\
\end{tabular}}
\caption{\textbf{Parameter table}, for both ODE and PDE models from Eqs. \ref{EQ_ODE_model} and \ref{EQ_PDE_model}.}
\label{table param}
\end{table}

Besides, for the initial conditions on the PDE model, we have set 
\begin{equation}
\begin{array}{l}
K(x)=B_{\mathrm{max}}\big(0.1\,\mathcal{G}(x;0.25,\sigma)+0.65\,\mathcal{G}(x;0.5,\sigma)+0.25\,\mathcal{G}(x;0.75,\sigma)\big),\vspace{5pt}\\

L(0,x)=L_0\big(0.5\,\mathcal{G}(x;0.4,\sigma)+0.5\,\mathcal{G}(x;0.6,\sigma)\big),\vspace{5pt}\\

B(0,x)=B_0\mathcal{G}(x;0,\sigma),
\end{array}
\end{equation}
with $L_0=10^7$ cells, $B_0=10^7$ cells, parameter $\sigma_i=\sigma_j=\sigma=1/20$ cells, without loss of generality, and where $\mathcal{G}(\cdot;\mu,\sigma)$ denotes a Gaussian kernel. The choice of the proportions of the maximal carrying capacity of $K(x)$ is chosen based on the Pre-B, Pro-B and transition proportions of B cells in healthy patients \cite{CHULIAN2021110685}. For the case of leukemic cells, we have chosen the possibility of two leukemic clones, one centered at normalized antigen expression of 0.4, and other at 0.6. We have, however, provided several simulations where we have included only one of such clones, to compare the behavior of the tumor depending on such expression levels, whereas the rest of simulations including one or two clones can be found in the Supplementary Information (Section \ref{SI}).

Both ODE- and PDE-like antigen-structured CAR T models from Eqs. \eqref{EQ_ODE_model} and \eqref{EQ_PDE_model} were implemented in MATLAB, where antigen expression was discretized for the PDE model and the resulting semi-discrete system is integrated in time. In the PDE model, the antigen domain was normalized to $x\in[0,1]$ and represented on a uniform grid with $N_x=100$ nodes and spacing $\Delta x=1/N_x$. Time integration was performed with MATLAB \texttt{ode45} over a daily output grid $t=0,1,\dots,t_{\mathrm{end}}$. The default horizon was $t_{\mathrm{end}}=300$ days, and interactive exploration was performed by varying selected parameters through MATLAB UI sliders. Global antigen burdens integral were implemented numerically by trapezoidal quadrature to optimize time computation.
To preserve biological interpretability and avoid nonphysical numerical artifacts, values below one cell unit in $L$ and $B$ were truncated to zero during integration. We have incorporated the concept of Minimal Residual Disease (MRD), which is usually studied in the bone marrow to detect leukemic cell levels \cite{campana2010minimal}. The positivity of such parameters is defined by the presence of more than $0.01\%$ leukemic cells in the bone marrow. For the sake of visualization and without loss of generality, we have chosen a threshold of $0
.1\%$ leukemic cells in the bone marrow for the subsequent simulations.

\section{Simulation of biological scenarios}

\subsection{ODE model simulations}
\label{Sec ODE model simulations}
Here we present specific scenarios to understand the general behavior of the models, specifically the ODE model presented in Section \ref{ODE model}. The model from Eqs. \eqref{EQ_ODE_model} specifically describes several cases. We omit the control case, where B cells are able to recover in absence of both CAR T cells and leukemic cells, as they behave normally as a logistic equation \cite{CHULIAN2021110685}. This results in 3 different cases to simulate: First cases 1-2 show how the model behaves similar to the ones present in the literature, where they only consider CD19$^+$ cells. A last case 3 shows how this behavior changes when in presence of a CD19$^-$ clone. We summarize the interpretation of all cases, shown in Figs.  \ref{fig:ODE_control} and \ref{fig:ODE_invasion}: 
\begin{itemize}

    \item[Case 1] \textbf{CAR T cells efficiently eliminate or control CD19$^+$ leukemic cells};

In this case, the ODE model resonates with other models as \cite{leon2021car} or \cite{martinezrubio2021}, as we set $L_N=0, \epsilon=0, \delta=0$. The $L_P$ population decreases at the time of CAR T infusion, and yet $C_A$ and $C_M$ cells begin to decline at the early stages. As a result, $L_P$ cells grow back significantly depending on the parameter $\alpha$, but in this model, also parameter $k$ and $h$ take relevance. 

\begin{itemize}
    \item 
For low $k$, this is, a lower activation threshold for the CAR T cells rate (see Figure \ref{fig:ODE_control} A.1 and A.2), cell population $L_P$ becomes extinct, where $B$ cells are also under attack of the $C_A$ cells. These last stabilize after the minimal expansion of $C_M$ cells, and are able to eliminate the CD19$^+$ tumor cells $L_P$. 

\item 
For higher $k$, such activation threshold allows for leukemic cells to grow more before CAR T cells expand, as they reach levels higher than the MRD level of the bone marrow. Here, parameter $h$, which is the treatment efficacy activation threshold, influences the dynamics of CAR T cells. A lower $h$ would imply that the treatment is able to be more efficient, even if they would need low levels of antigen to promote the killing of $L_p$ cells. In this sense, Figure \ref{fig:ODE_control} B.1 and B.2 show that CAR T cells oscillate more, but create lower peaks of CAR T needed to attack leukemic cells. This results in leukemic cells having a higher amount of oscillations. On the other hand, Figure \ref{fig:ODE_control} C.1 and C.2 show that a higher $h$ regulates CAR T cells to have less oscillations and stabilize in higher levels, and thus killing efficacy is lower, needing more activated CAR T cells. These are both plausible lotka-volterra dynamics where $L_P$ cells are the controlled targets of $C_A$ cells, and $C_M$ are the memory cells that can maintain the amount of T cells present in the body.

This results that a lower $k$ would effectively help as an activation threshold for $C_A$ to replicate and thus deplete antigen positive leukemic cells $L_P$. This would then imply a lower amount of healthy $B$ cells, and the patient would be considered as immunodeficient.  
\end{itemize}

    \item[Case 2] \textbf{CAR T cells are not able to reduce CD19$^+$ leukemic cells};

For a high $k$ and $\rho_{C_A}$, CAR T cells are not able to activate themselves enough in the presence of a high amount of either $B$ or $L_P$ cells. In this sense, CAR T cells $C_A$ would not replicate enough to counterattack the presence of leukemic antigen positive cells. This can be observed in Figure \ref{fig:ODE_invasion} A.1 and A.2, where CAR T cells cannot expand enough. Even if they do, when leukemic cells grow higher, they are not enough to control the disease at reasonable levels. 
     \item[Case 3]  \textbf{CAR T cells are not able to control the rapid growth of CD19$^-$ cells}

The inclusion of $L_N$ cells (CD19$^-$  antigen expressed cells) changes the behavior of the ODE model. In this sense, as seen in Figure \ref{fig:ODE_invasion} B.1 and B.2,  the minimal amount of possible $L_N$ cells (by considering $L_N(0)=1$ cell and setting parameter $\delta>0$) would imply an exponential growth of the CD19$^-$ cells $L_N$ at a certain moment. Parameter $\epsilon$ is directly related on when the CD19$^-$ recurrence will happen, as the higher the $\epsilon$ value, the later in time the  CD19$^-$ relapse.

\begin{figure}[H]
    \centering
    \includegraphics[width=\linewidth]{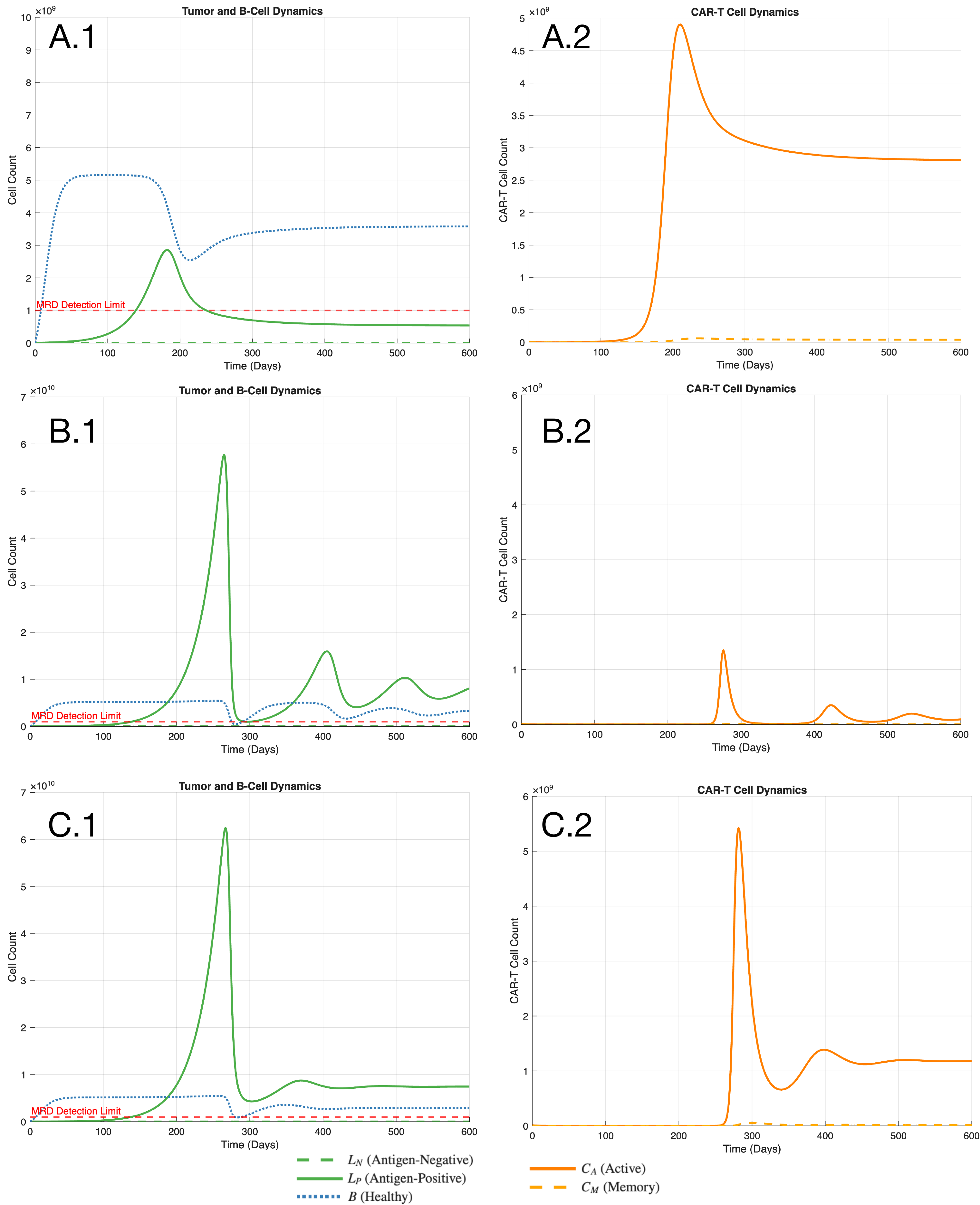}
    \caption{\textbf{Case 1. Elimination or control of leukemic CD19$^+$ cells in ODE Model from Eqs. \eqref{EQ_ODE_model}}. Left panel: leukemic CD19$^+$ cells $L_P$ (green solid line), healthy $B$ cells (blue dotted lines), CD19$^-$ $L_N$ set to 0 (green, dashed line) and threshold bone marrow level (red horizontal dashed line). Right panel: CAR T activated cells $C_A$ (orange solid lines) and memory $C_M$ cells (orange dashed lines). Each row represent different parameter sets: A with $k=2\cdot 10^{10}$ cells and $h=10^{11} $ cells; for $k=5\cdot10^{10} $ cells, both rows B and C respectively show $h=10^{9} $ cells or $h=10^{11} $ cells.}
    \label{fig:ODE_control}
\end{figure}

\begin{figure}[H]
    \centering
    \includegraphics[width=\linewidth]{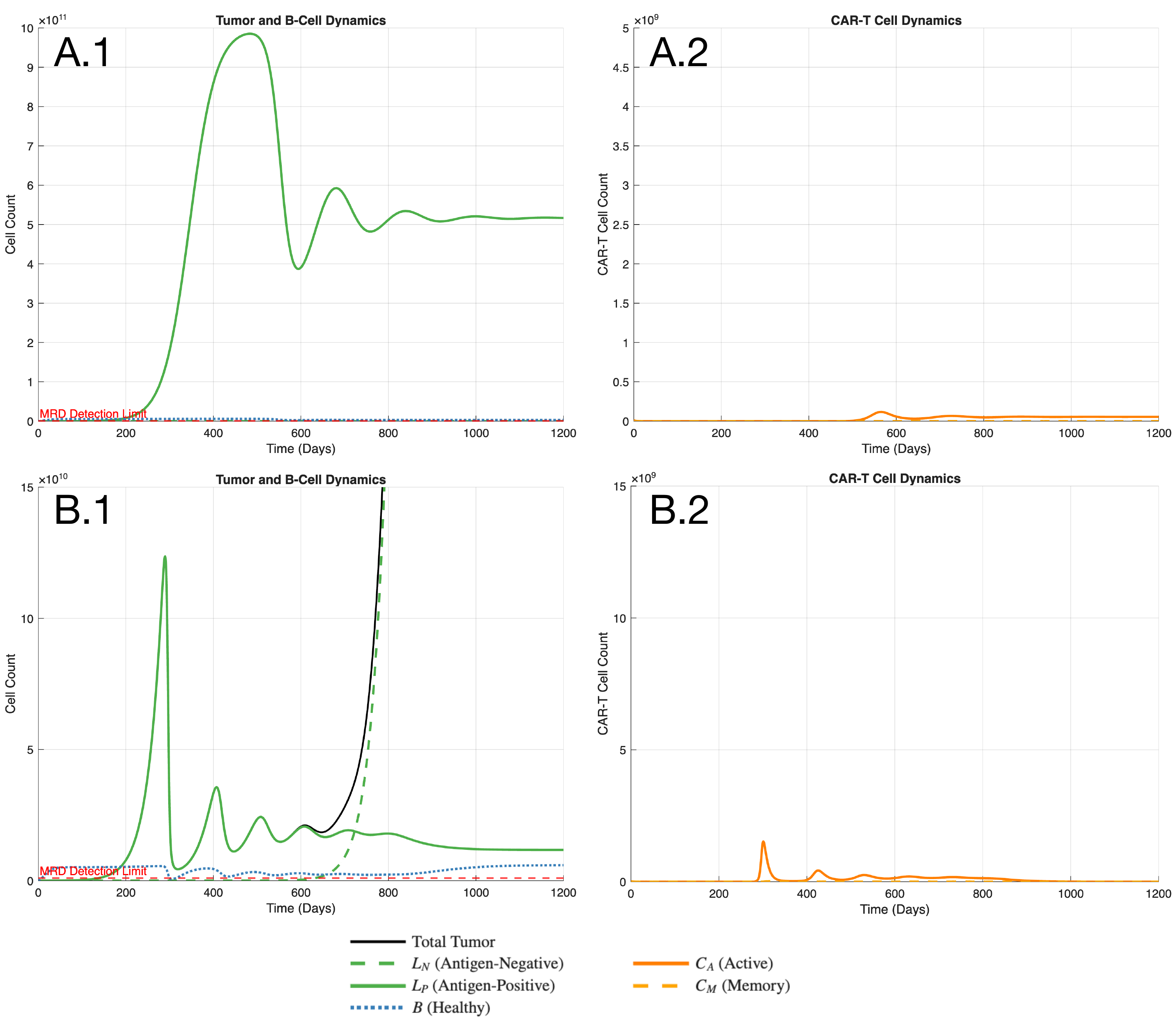}
    \caption{\textbf{Uncontrolled growth of tumor cells (Cases 2 and 3) in ODE Model from Eqs. \eqref{EQ_ODE_model}.} In row A, case 2 shows how leukemic CD19$^+$ cells grow due to a lower reproduction rate of activated CAR T cells $\rho_{C_A}=0.9/2 \ day^{-1}$ and a higher activation threshold level $k=10^{12} $ cells. In row B, we set $L_N(0)=1 \ cell $, this is, just one CD19$^-$ cells. CD19$^+$ tumor levels are controlled for some time, but then an exponential growth of $L_N$ cells outgrows the $L_P$ cells, and CAR T cells are lost. Left panel: total tumor (black solid thin line), leukemic CD19$^+$ cells $L_P$ (green solid thick line), healthy $B$ cells (blue dotted lines), CD19$^-$ $L_N$ (green, dashed line) and MRD  bone marrow threshold level (red horizontal dashed line). Right panel: CAR T activated cells $C_A$ (orange solid lines) and memory $C_M$ cells (orange dashed lines). All  have $h=10^9 $ cells and other parameters are taken from Table \ref{table param}.}
    \label{fig:ODE_invasion}
\end{figure}

\end{itemize}

\subsection{PDE model simulations}

This section shows PDE simulations designed to characterize baseline hematopoietic dynamics and relapse patterns under different immune-pressure regimes from Eqs. \eqref{EQ_PDE_model}.

\subsubsection*{Baseline dynamics: recovery and leukemia invasion}
\begin{figure}[H]
    \centering
    \includegraphics[width=0.45\textwidth]{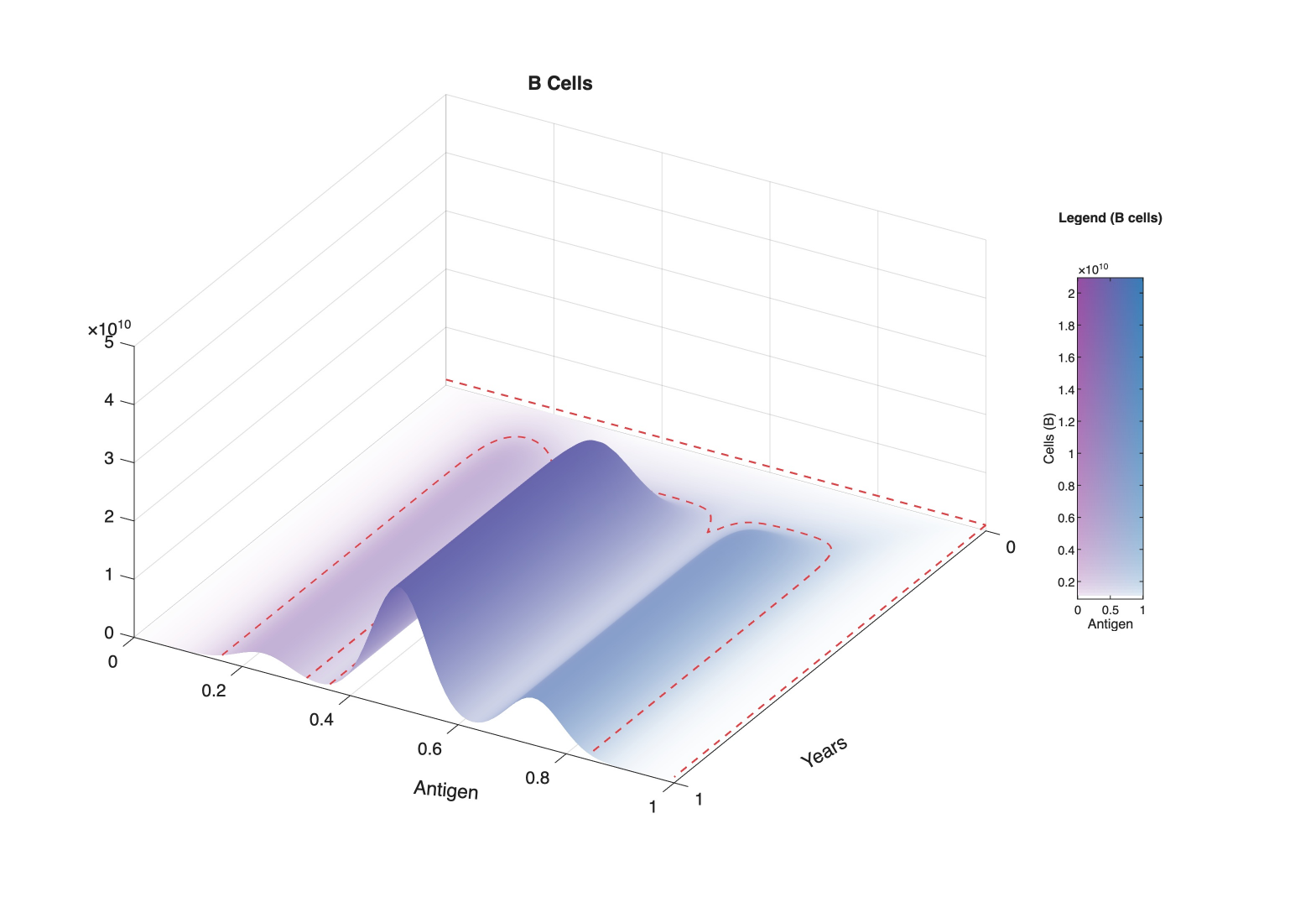} \includegraphics[width=0.45\textwidth]{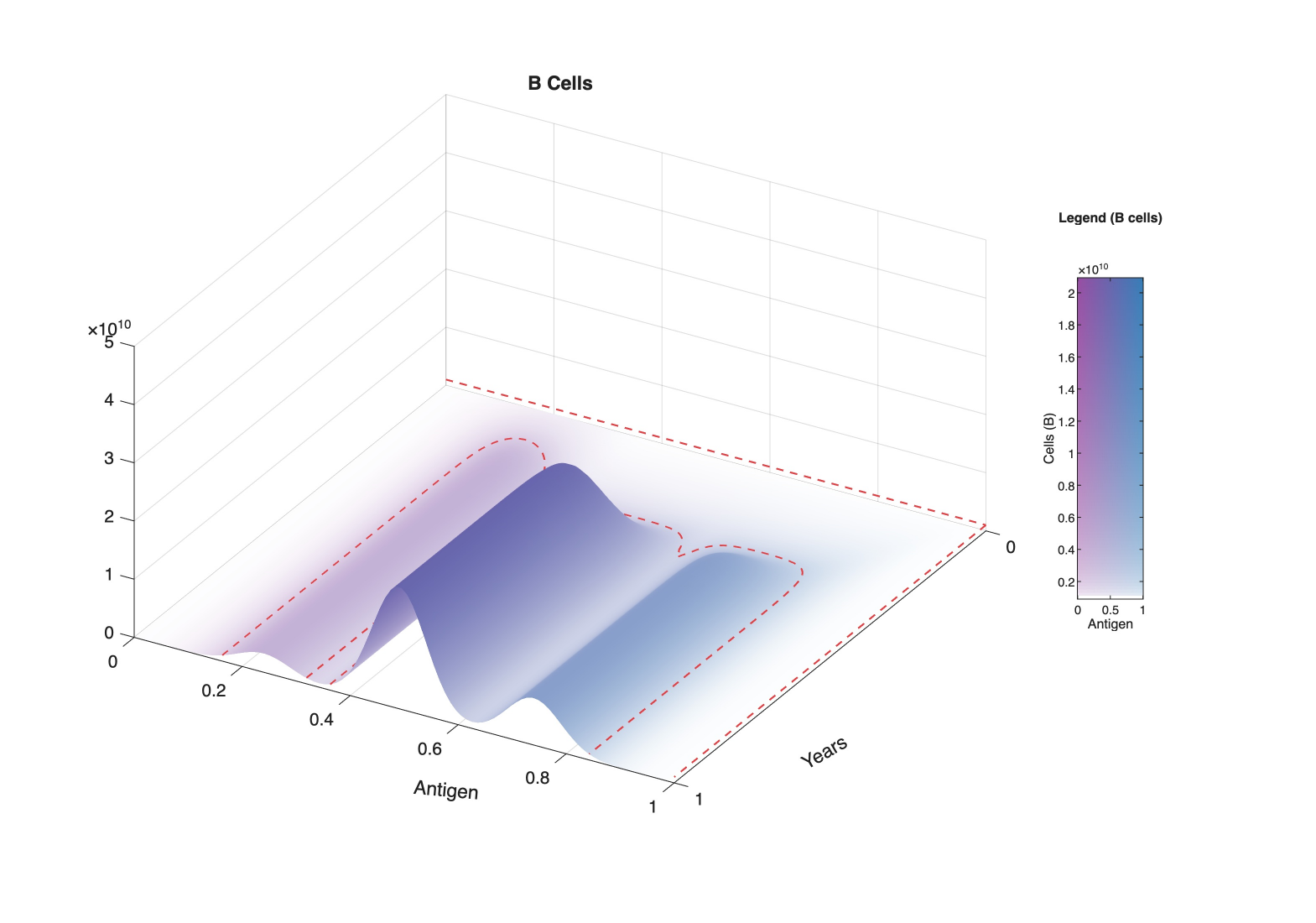} 
    \caption{\textbf{Normal hematopoiesis simulations in PDE Model from Eqs. \eqref{EQ_PDE_model}.} We use the parameters $s=0.01$, $s=0.001$, respectively.  $L_0=C_{A0}=C_{M0}=0 $ cells , and simulate it for 1 year. Red dashed lines represents the plane where cells achieve 0.1\% of the total bone marrow count.}
    \label{fig:Scenario0_300}
\end{figure}

 We first analyze control configurations without CAR T cells so that we can establish reference behavior. In the physiological hematopoiesis scenario, the initial distribution $B(0,x)$ evolves toward homeostasis. The parameter $s$ controls the time scale of recovery: larger $s$ yields faster convergence to steady B cell levels (Figure \ref{fig:Scenario0_300}).

\begin{figure}[H]
    \centering
    \includegraphics[width=0.45\textwidth]{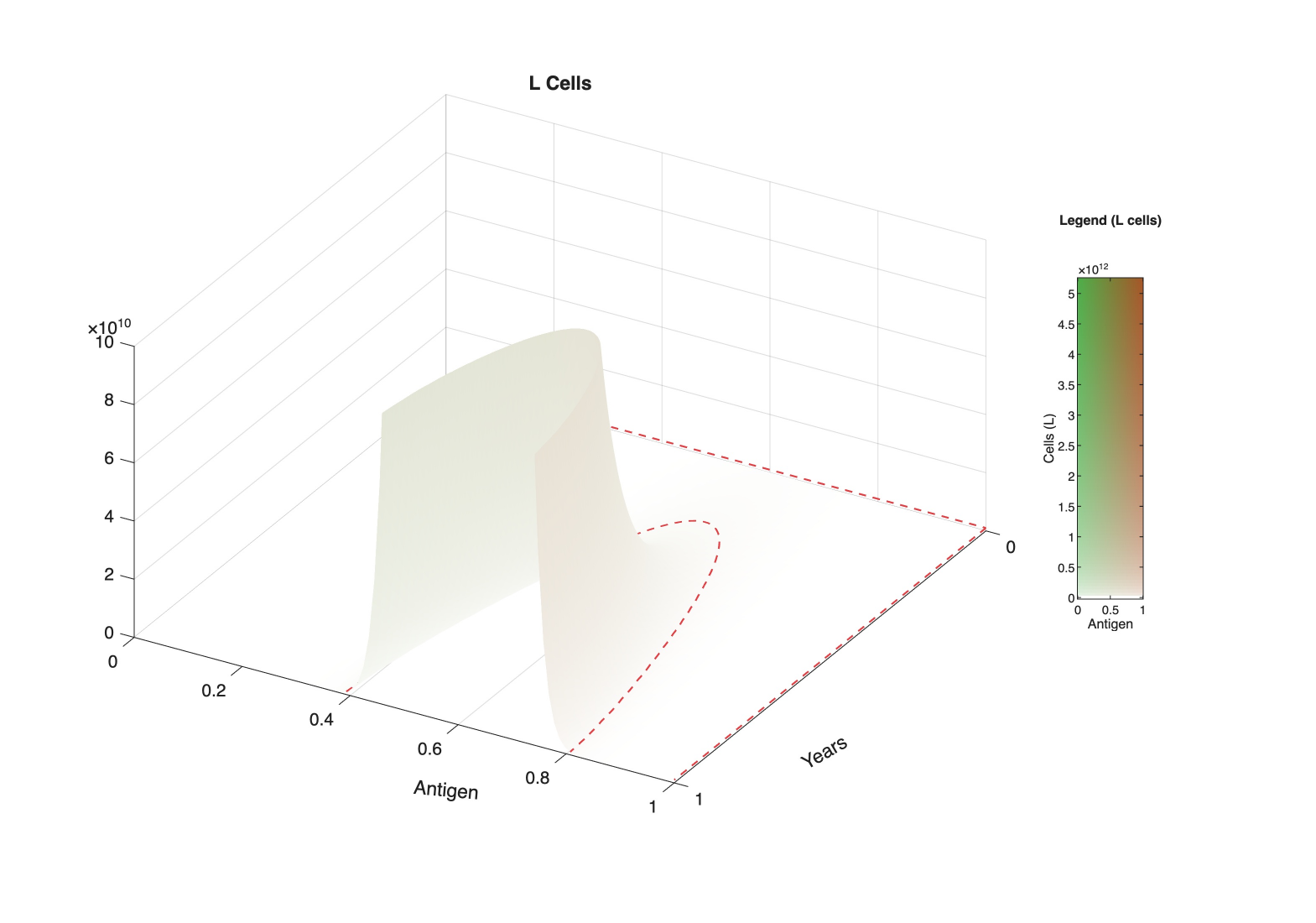}
    \includegraphics[width=0.45\textwidth]{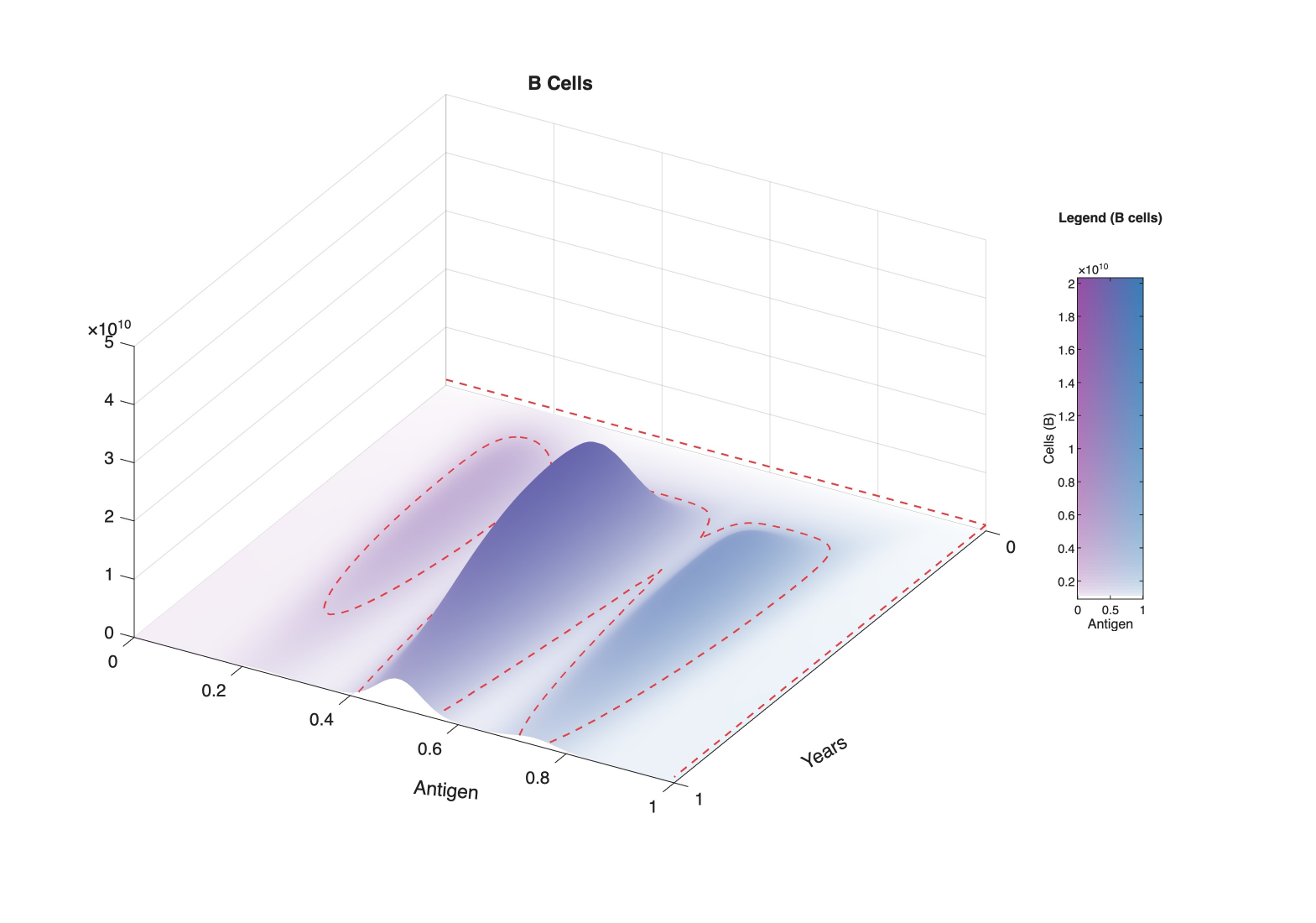}
    \caption{\textbf{Leukemic invasion and B cells decay with a one heterogeneous clone in PDE Model from Eqs. \eqref{EQ_PDE_model}.} We use parameter set $s=0.01$, $L_0=B_0=10^7$ cells, $ C_{A0}= C_{M0}=0$ cells, and simulate it for 1 year. Red dashed lines represents the plane where cells achieve 0.1\% of the total bone marrow count.}
    \label{fig:S0_leukemia_1_clon}
\end{figure}

We then evaluate leukemic invasion in the absence of CAR T treatment. For a single leukemic clone centered at $x_1=0.6$, total leukemic burden increases around day 100, accompanied by progressive depletion of healthy B cells (Figure \ref{fig:S0_leukemia_1_clon}). The same qualitative behavior is preserved under a more heterogeneous initial leukemic distribution (see GUIs from Supplementary Information, Section \ref{SI}).

\subsubsection*{CD19$^+$ relapse under CAR T pressure}

For the CD19$^+$ relapse regime, we consider $h=0 $, a well as a single Gaussian initial leukemic profile $L(0,x)$, and a high activation-threshold setting $k=10^{10}$ cells . We compare two scenarios depending on the initial disease burden $L_0=10^7 $ or $10^9 $ cells.

\begin{figure}[H]
\vspace{-0.4cm}
    \centering
    
    \includegraphics[width=0.95\linewidth]
    {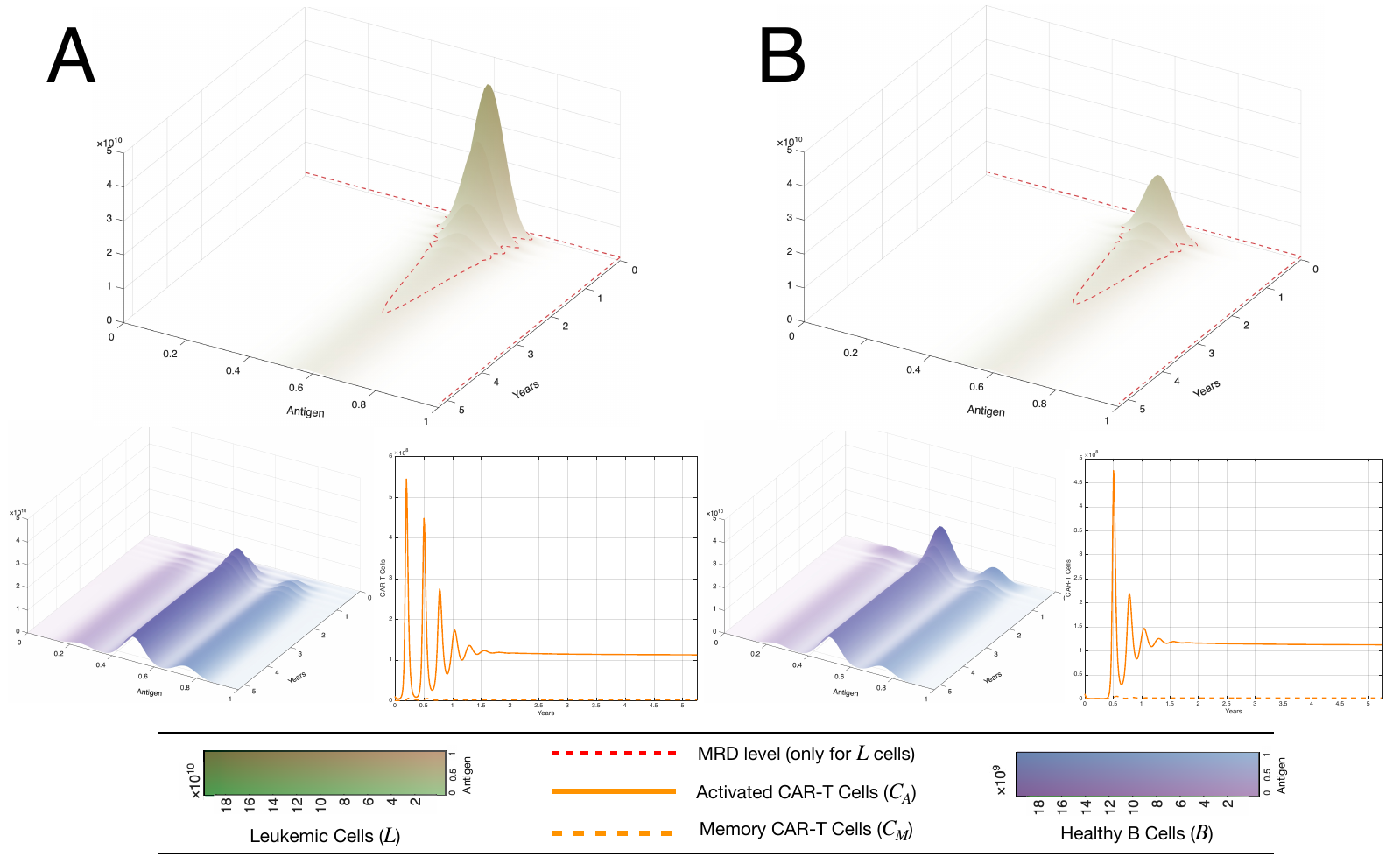}
    \caption{\textbf{PDE model  dynamics from Eqs. \eqref{EQ_PDE_model} under a CD19$^+$ leukemia.} We use model from \eqref{PDE model} and simulate two different scenarios \textbf{(A)} with high $L_0=10^9 $ cells , and \textbf{(B)} with $L_0=10^7 $ cells . Both simulations have a single clone centered at $x_0=0.6$ for the antigen expression, as well as the parameters used in Table \ref{table param} specifically for $h=0$.}
    \label{fig:CD19pos}
\end{figure}

\begin{itemize}
\vspace{-5pt}
    \item For higher $L_0$ (see Figure \ref{fig:CD19pos}A), CAR T activation is increased, occurs in sharper bursts, producing larger amplitude oscillations in activated CAR T cells. This leads to very oscillating lotka-volterra dynamics, present between activated CAR T cells and both leukemic and B cells. In this scenario, healthy B cells expansion is lowered by the rapid growth of CAR T cells, but it stabilizes over time.

    \item For lower $L_0$ (see Figure \ref{fig:CD19pos}B), the first expansion of CAR T cells is not as high and quick in comparison to a higher amount of initial leukemic cells $L_0$. In fact, B cells are able to reactivate themselves prior to the leukemic expansion, which is lower in comparison to a high $L_0$. CAR T cells expand later in time and deplete more smoothly the leukemic burden, reaching also an stabilization of the B cell levels.

\end{itemize}

These qualitative trends persist for lower high activation-threshold setting  $k=10^{9}$ cells and also for the presence of two distinct leukemic clones, one centered at $x_0=0.4$ and other centered at $x_0=0.6$. Such dynamics can be observed in the Supplementary Information (Section \ref{SI}).

\subsubsection*{CD19$^-$ relapse and immune escape}

\begin{figure}[H]
    \centering
    \includegraphics[width=1\linewidth]%
    {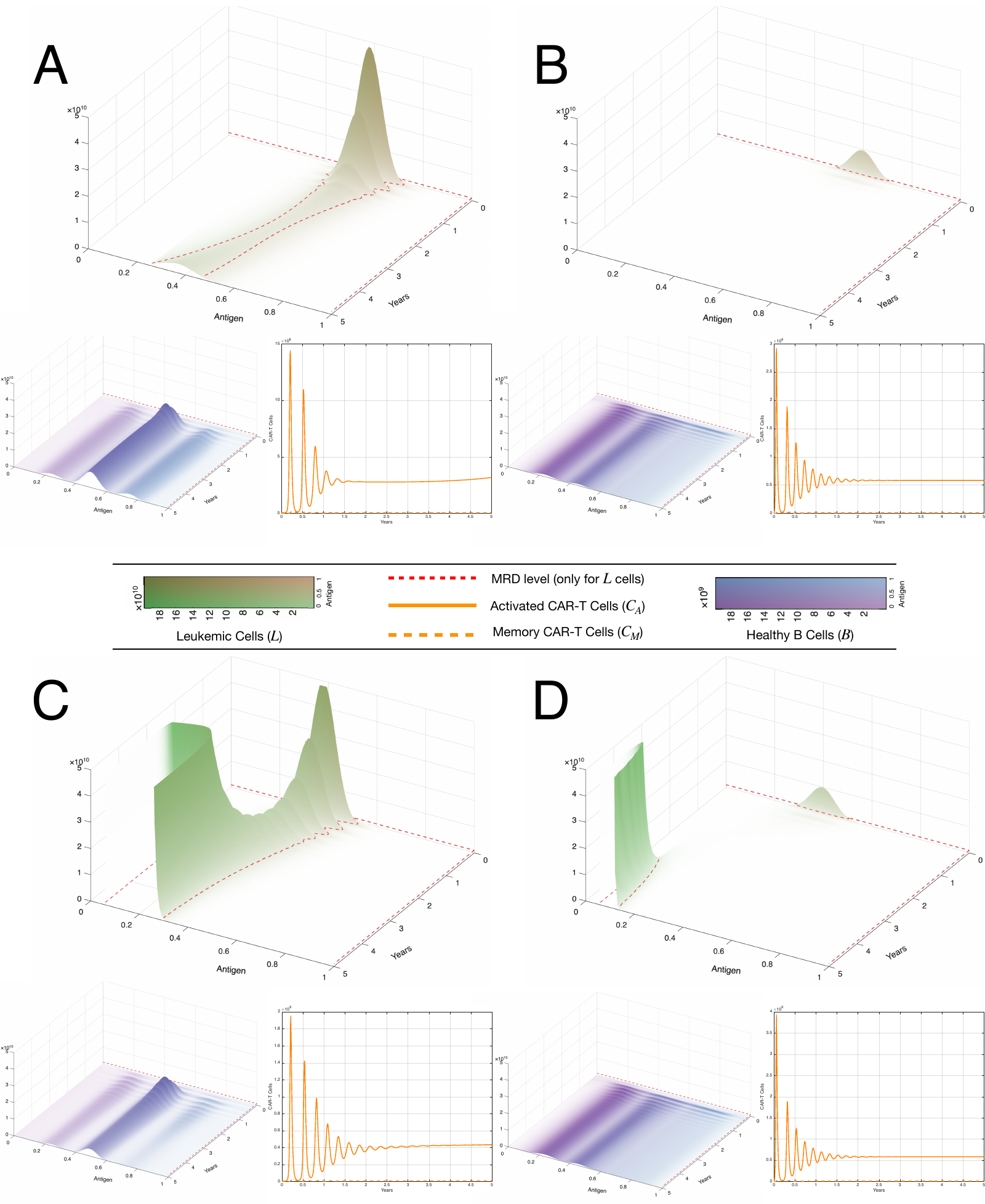}
    \caption{\textbf{PDE model  dynamics from Eqs. \eqref{EQ_PDE_model} under the possibility of CD19$^-$ leukemia.} We set specifically $h=0.75$ to show the possibility of antigen immune escape.  include different possibilities for parameters $x_0$ (where the leukemic clone is centered) and CAR T activation threshold $k$:  \textbf{(A)} $x_0=0.6, k=10^{10}$. \textbf{(B)} $x_0=0.6, k=10^{9}$. 
    \textbf{(C)} $x_0=0.4, k=10^{10}$. 
    \textbf{(D)} $x_0=0.4, k=10^{9}$. All other parameters are the ones used in Table \ref{table param}.}
    \label{fig:PDE_scenarios}
\end{figure}

To study CD19$^-$ relapse and immune escape, we set $h=0.75$ and analyze high-burden initial conditions ($L_0=10^9 $ cells) under heterogeneous antigen-distribution structures, and are shown in Figure \ref{fig:PDE_scenarios}.

\begin{itemize}
    \item With a single high-antigen leukemic clone (Gaussian centered at $x_0=0.6$),  CAR T cells effectively recognize the leukemic population, but it depends also on higher or lower $k$ values:
    \begin{itemize}
        \item For higher $k$ values (for example, $k=10^{10} $ cells, see Figure \ref{fig:PDE_scenarios}A, CAR T cells activate slower but with very high peaks, as they need higher values of antigen to reactivate. In fact, CAR T cells are able to control the disease, but several relapses occur along time, slowly selecting leukemic cells to a lower antigen expression. B cells are able to grow back, but at cost of leukemia re-expansion.

        \item For lower $k$ values (for example, $k=10^9$ cells, see Figure \ref{fig:PDE_scenarios}B), however, even if leukemia is able to grow back, CAR T cells are able to reactivate themselves at lower values, creating lower, more oscillating peaks of activated CAR T cells. These CAR T cells, even if at lower stabilization levels, control better the leukemia burden, and the corresponding B cells are then depleted due to their presence.
    \end{itemize}

    \item With single low-antigen leukemic clone (Gaussian centered at $x_0=0.4$) outcomes also become $k$-dependent, and in both scenarios CAR T cells are not able to control the growth of leukemic CD19$^{-}$ cells:
    \begin{itemize}
        \item With higher values of $k$ ($k=10^{10}$ cells, see Figure \ref{fig:PDE_scenarios}C) CD19$^-$ relapses are way more aggressive, with higher peaks of growth and an uncontrolled behavior, sooner in comparison to a clone centered at higher values (for example, as shown previously with $x_0=0.6$). The recovery of healthy B cells from the bone marrow follows. 
        \item With lower values of $k$ ($k=10^9$ cells, see Figure \ref{fig:PDE_scenarios}D), CAR T cells successfully eradicates the leukemic population during the very first years. However, the repopulation of leukemic cells are biased toward lower-CD19 expression, and this selection shifts the leukemic distribution toward low-antigen states and eventually yields CD19$^-$ escape. This late relapse is preceded by CAR T contraction and transient healthy B cell recovery, indicating insufficient CAR T reactivation.
    \end{itemize}
    
\end{itemize}

The same mechanisms are observed for lower initial leukemic burden (for example with $L_0=10^7$ cells, see Supplementary Information, Section \ref{SI}). 

\section{Sensitivity analyses}

To quantify parameter influence over time, we performed a variance-based global sensitivity analysis using Sobol indices with Saltelli sampling. The uncertain parameter vector for the ODE model was defined as:
$$\theta_{ODE}=(B_{\text{Stem}}, B_{\text{max}}, k, \rho_{CA}, \gamma_{AM}, \tau_{CA}, \alpha, h, B_0, L_0, \epsilon, \delta),$$
while for the PDE model, the vector included spatial distribution parameters:
$$\theta_{PDE}=(B_{\text{Stem}}, B_{\text{max}}, s, x_0, m, k, \rho_{CA}, \gamma_{AM}, \tau_{CA}, \alpha, h, B_0, L_0).$$
Parameters such as $k$, $h$, and $L_0$ were sampled on a logarithmic scale to explore several orders of magnitude, reflecting the high uncertainty in biological rates and initial tumor burden.

We generated two independent sample matrices $A,B\in\mathbb{R}^{N\times d}$ with $N=1000$ and $d=14$ (ODE) or $d=15$ (PDE), and constructed hybrid matrices $A_B^{(i)}$ by replacing the $i$-th column of $B$ with the corresponding column of $A$. Therefore, the total number of model evaluations was $N_{\mathrm{eval}}=N(2+d)$. Each parameter set was simulated using \texttt{ode15s} (\texttt{RelTol}$=10^{-5}$, \texttt{AbsTol}$=10^{-7}$) with non-negativity constraints. For the PDE model, the antigen expression grid was discretized into $N_x=50$ points. The outputs analyzed were: activated CAR T cells ($C_A$), memory CAR T cells ($C_M$), total leukemic load ($L$), total healthy B cells ($B$), Antigen-Negative subsets (with $x<0.5$ for PDE and $L_N$ for ODE), and  Antigen-Positive subsets (with $x\geq0.5$ for PDE and $L_P$ for ODE).

For each output and time point, first-order and total-order Sobol indices were computed as
\[
S_i(t)=\frac{\mathbb{E}\left[Y_B(t)\left(Y_{A_B^{(i)}}(t)-Y_A(t)\right)\right]}{\mathrm{Var}(Y(t))},
\]
\[
S{T_i}(t)=\frac{\mathbb{E}\left[\left(Y_A(t)-Y_{A_B^{(i)}}(t)\right)^2\right]}{2\,\mathrm{Var}(Y(t))}.
\]
Following standard numerical practice, small negative estimates due to Monte Carlo noise were truncated to zero. Indices were presented as time-resolved stacked profiles to compare transient versus long-time influence and to separate main effects ($S_i$) from interaction-driven effects ($S{T_i}-S_i$).

Both the PDE and ODE models share several consistent parameter ranges, including:
$B_{\text{Stem}}, B_0 \in [10^6, 10^9]$, $B_{\text{max}} \in [10^9, 10^{11}]$, $k \in [10^8, 10^{12}]$, 
$\rho_{CA} \in [0.01, 1]$, $\gamma_{AM} \in [10^{-4}, 10^{-2}]$, $\tau_{CA} \in [1, 15]$, 
$L_0 \in [10^3, 10^9]$, and $\alpha \in [3 \times 10^{-11}, 3 \times 10^{-9}]$. 
However, they diverge in the definition of $h$, which ranges from $[0, 1]$ in the PDE model 
and $[10^9, 10^{12}]$ in the ODE model. Additionally, the PDE model uniquely includes 
spatial parameters $s \in [10^{-4}, 10^{-1}]$, $x_0 \in [0, 100]$, and $m \in [0, 10]$, 
while the ODE model exclusively incorporates $\epsilon \in [0, 10^{-2}]$ and $\delta \in [0, 10^{-15}]$.

\begin{figure}[H]
     \centering

\includegraphics[width=\textwidth]{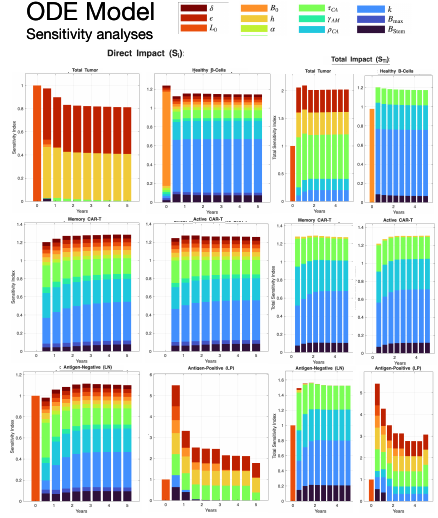}
    \caption{\textbf{First and Total order sensitivity analysis (ODE model from \eqref{EQ_ODE_model})}.  The graphics show the first-order $S_i$ (direct impact) total order $ST_i$ (total impact) indices over different years, from 0 to 5, on the different subpopulations for the model, $L=L_N+L_P,B,C_M,C_A,L_N$ and $L_P$. The parameters observed were $B_{\text{Stem}}, B_0, B_{\text{max}}, k, \rho_{CA}, \gamma_{AM} ,\tau_{CA} , L_0 $ and $\alpha $, and specifically for the ODE $\epsilon,\delta$ and $h$, as this last has a different range for the ODE model.}
    \label{fig:ODE_Rainbow}
\end{figure}

\begin{figure}[H]
    \centering
    
\includegraphics[width=\textwidth]{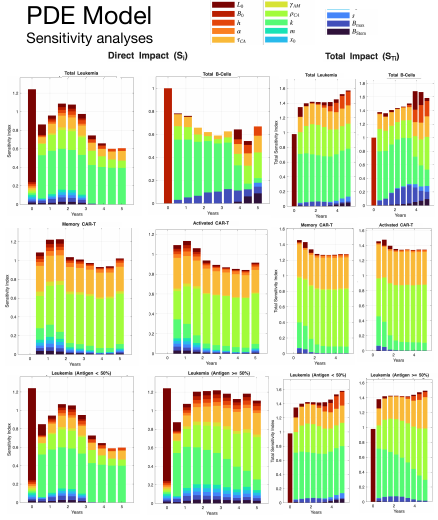}
    \caption{\textbf{First and Total order sensitivity analysis (PDE model from \eqref{EQ_PDE_model})}. The graphics show the first-order $S_i$ (direct impact) total order $ST_i$ (total impact) indices over different years, from 0 to 5, on the different subpopulations for the model, $L,B,C_M,C_A$, and includes the influence specifically for antigen negative leukemic cells, with $x<0.5$, and for antigen positive leukemic cells, with $x\geq0.5$. The parameters observed were $B_{\text{Stem}}, B_0, B_{\text{max}}, k, \rho_{CA}, \gamma_{AM} ,\tau_{CA} , L_0 $ and $\alpha $, and specifically for the PDE $s,x_0,m$  and $h$, as this last has a different range for the PDE model. }
    \label{fig:PDE_Rainbow}
\end{figure}

We have included the sensitivity analyses results for both the ODE model in Figure \ref{fig:ODE_Rainbow}, and for the PDE model in Figure \ref{fig:PDE_Rainbow}. For the sake of readability, we have only included the results in heatmaps for total indices $S{T_i}$ in the Supplementary Information (Section \ref{SI}). 
We also present a summary of such results in Table \ref{tab:sensitivity_full}, where we show the three most relevant parameters over time for each subpopulation for the total indices $S{T_i}$.

\begin{table}[htbp]
\centering
\scriptsize
\begin{tabular}{l | ccc | ccc}
\hline
 & \multicolumn{3}{c|}{\textbf{ODE MODEL}} & \multicolumn{3}{c}{\textbf{PDE MODEL}} \\
\textbf{Timeline} & \textbf{Rank 1} & \textbf{Rank 2} & \textbf{Rank 3} & \textbf{Rank 1} & \textbf{Rank 2} & \textbf{Rank 3} \\ 
\hline
\multicolumn{7}{c}{\textit{Target: Total Leukemia / Total Tumor}} \\
\hline
Year 1 & $\tau_{CA}$ (0.84) & $h$ (0.44) & $\epsilon$ (0.44) & $k$ (0.66) & $\rho_{CA}$ (0.37) & $\tau_{CA}$ (0.27) \\
Year 3 & $\tau_{CA}$ (0.80) & $h$ (0.41) & $\epsilon$ (0.41) & $k$ (0.62) & $\rho_{CA}$ (0.40) & $\tau_{CA}$ (0.26) \\
Year 5 & $\tau_{CA}$ (0.80) & $h$ (0.41) & $\epsilon$ (0.41) & $k$ (0.65) & $\rho_{CA}$ (0.34) & $\tau_{CA}$ (0.25) \\
\hline
\multicolumn{7}{c}{\textit{Target: B cells}} \\
\hline
Year 1 & $k$ (0.69) & $\rho_{CA}$ (0.26) & $\tau_{CA}$ (0.16) & $k$ (0.70) & $\rho_{CA}$ (0.30) & $\tau_{CA}$ (0.20) \\
Year 3 & $k$ (0.69) & $\rho_{CA}$ (0.26) & $\tau_{CA}$ (0.16) & $k$ (0.71) & $\rho_{CA}$ (0.26) & $B_{\text{max}}$ (0.25) \\
Year 5 & $k$ (0.69) & $\rho_{CA}$ (0.26) & $\tau_{CA}$ (0.16) & $\rho_{CA}$ (0.48) & $\tau_{CA}$ (0.30) & $k$ (0.24) \\
\hline
\multicolumn{7}{c}{\textit{Target: Activated CAR T ($C_A$)}} \\
\hline
Year 1 & $k$ (0.54) & $\rho_{CA}$ (0.35) & $\tau_{CA}$ (0.27) & $\rho_{CA}$ (0.59) & $\tau_{CA}$ (0.40) & $k$ (0.33) \\
Year 3 & $k$ (0.59) & $\rho_{CA}$ (0.34) & $\tau_{CA}$ (0.25) & $\rho_{CA}$ (0.75) & $\tau_{CA}$ (0.44) & $k$ (0.11) \\
Year 5 & $k$ (0.59) & $\rho_{CA}$ (0.34) & $\tau_{CA}$ (0.24) & $\rho_{CA}$ (0.77) & $\tau_{CA}$ (0.44) & $k$ (0.10) \\
\hline
\multicolumn{7}{c}{\textit{Target: Memory CAR T ($C_M$)}} \\
\hline
Year 1 & $k$ (0.50) & $\rho_{CA}$ (0.38) & $\tau_{CA}$ (0.27) & $\rho_{CA}$ (0.59) & $\tau_{CA}$ (0.37) & $k$ (0.30) \\
Year 3 & $k$ (0.57) & $\rho_{CA}$ (0.34) & $\tau_{CA}$ (0.24) & $\rho_{CA}$ (0.72) & $\tau_{CA}$ (0.40) & $k$ (0.09) \\
Year 5 & $k$ (0.57) & $\rho_{CA}$ (0.33) & $\tau_{CA}$ (0.24) & $\rho_{CA}$ (0.74) & $\tau_{CA}$ (0.40) & $k$ (0.08) \\
\hline
\multicolumn{7}{c}{\textit{Target: Antigen-Negative (ODE: $L_N$ | PDE: Antigen < 50\%)}} \\
\hline
Year 1 & $k$ (0.51) & $\rho_{CA}$ (0.43) & $\tau_{CA}$ (0.38) & $k$ (0.67) & $\rho_{CA}$ (0.37) & $\tau_{CA}$ (0.27) \\
Year 3 & $k$ (0.59) & $\rho_{CA}$ (0.40) & $\tau_{CA}$ (0.32) & $k$ (0.63) & $\rho_{CA}$ (0.39) & $\tau_{CA}$ (0.26) \\
Year 5 & $k$ (0.59) & $\rho_{CA}$ (0.41) & $\tau_{CA}$ (0.31) & $k$ (0.66) & $\rho_{CA}$ (0.33) & $\tau_{CA}$ (0.24) \\
\hline
\multicolumn{7}{c}{\textit{Target: Antigen-Positive (ODE: $L_P$ | PDE: Antigen $\geq$ 50\%)}} \\
\hline
Year 1 & $\tau_{CA}$ (0.80) & $h$ (0.80) & $\epsilon$ (0.80) & $k$ (0.64) & $\rho_{CA}$ (0.40) & $\tau_{CA}$ (0.28) \\
Year 3 & $\tau_{CA}$ (0.70) & $h$ (0.70) & $\epsilon$ (0.70) & $k$ (0.53) & $\rho_{CA}$ (0.50) & $\tau_{CA}$ (0.32) \\
Year 5 & $\tau_{CA}$ (1.00) & $h$ (0.70) & $\epsilon$ (0.70) & $\rho_{CA}$ (0.69) & $\tau_{CA}$ (0.42) & $k$ (0.28) \\
\hline
\end{tabular}
\caption{\textbf{Comparative Sensitivity Analysis: Top 3 Parameters in ODE and PDE Models over time}. The results are taken from the sensitivity analyses from  \ref{fig:ODE_Rainbow} and \ref{fig:PDE_Rainbow} and only for the total indices $S{T_i}$.}
\label{tab:sensitivity_full}
\end{table}

The first- and total-order Sobol analyses reveal specifically for the PDE model a consistent behavior in parameter influence over time, though with critical shifts in the hierarchy of dominance compared to the ODE version. In the early phase, all four outputs show pronounced variability in sensitivity, reflecting strong nonlinear transients and competition among mechanisms as the CAR T population first encounters the distributed antigen profiles of the leukemia cells. As time progresses, the sensitivity structure becomes more regular, especially for total leukemic and total healthy B cell loads. This evolution indicates a transition from a system governed by initial conditions to one defined by a robust long-time influence hierarchy centered on the activation threshold and therapeutic persistence. Furthermore, the analyses confirm that total-order indices $S{T_i}$ are systematically larger than first-order indices $S_{i}$, demonstrating that interaction effects are substantial and treatment efficacy is a product of coupled mechanisms rather than isolated parameters.

The ODE model sensitivity analyses shown in Figure \ref{fig:ODE_Rainbow} illustrate a sensitivity profile heavily dominated by CAR T cell persistence and activation efficiency. For the total tumor population ($L_N + L_P$), the lifespan of activated CAR T cells $\tau_{CA}$ is the primary driver, maintaining a consistently high total-order index ($S_{Ti} \approx 0.80$) throughout the five-year simulation. This dominance is mirrored in the antigen-positive leukemia subpopulation $L_P$, where $\tau_{CA}$ reaches its maximum influence ($S_{Ti} = 1.00$) by Year 5, alongside significant contributions from the half-maximal activation constant $h$ and the small-scale constant $\epsilon$. Conversely, the healthy B cell population and the CAR T compartments (active and memory) are most sensitive to the activation threshold $k$, which maintains stable indices near 0.69 and 0.57, respectively. This suggests that in the ODE framework, the therapeutic outcome is strictly constrained by the discrete transition of cells into an active state and the subsequent longevity of that active population.

In contrast, the PDE model sensitivity analyses shown in Figure \ref{fig:PDE_Rainbow} reveal a shift toward expansion kinetics and activation thresholds as the primary governing mechanisms in a continuous antigen environment. The total leukemia load and its subpopulations ($<50\%$ and $\geq50\%$ antigen expression) are consistently dominated by the activation threshold $k$, which remains the highest-ranked parameter with indices between 0.62 and 0.67. This indicates that the bottleneck for tumor control in the PDE model is the collective stimulus required to trigger the CAR T response across the antigen spectrum. For the CAR T populations themselves, the expansion rate $\rho_{CA}$ is the overwhelming driver, with its sensitivity index increasing over time to reach 0.77 for active cells and 0.74 for memory cells by Year 5. This long-term evolution signifies that the PDE system's stability and its ability to prevent late-stage recurrence are more dependent on the proliferative capacity of the CAR T cells than on their individual lifespan or killing efficiency.

Ultimately, comparing the two models demonstrates that the continuous antigen spectrum in the PDE model shifts the therapeutic bottleneck. In the ODE model, the system is highly sensitive to the discrete persistence and efficiency of the CAR T cells ($\tau_{CA}$ and $h$). However, the PDE model highlights how system stability depends on the collective stimulus provided by the biomass, with $k$ remaining the master regulator of the tumor load. The emergence of $\rho_{CA}$ and $\tau_{CA}$ as vital drivers in the long-term PDE dynamics suggests that a successful clinical outcome requires a balance between a low activation threshold to detect low-burden disease and high expansion/persistence to prevent late-stage recurrence.

\section{Parameter influence on antigen-dependent relapse}

In this final section, we evaluate the long-term (10-year) dynamics of the PDE model by quantifying the total tumor burden and the evolution of mean antigen expression. These metrics are analyzed as a function of two critical parameters: $k$, representing the activation threshold of CAR T cells, and $h$, which dictates the selection pressure and the resulting phenotypic shift from antigen-positive to antigen-negative leukemic populations. While comprehensive analyses were performed for both the ODE and PDE models across several initial antigen distributions (see Supplementary Information, Section \ref{SI}), for clarity, we focus here on the PDE model results for a single leukemic clone initially centered at $x_0 = 0.6$.

\begin{figure}[H]
    \centering
    \includegraphics[width=\linewidth]{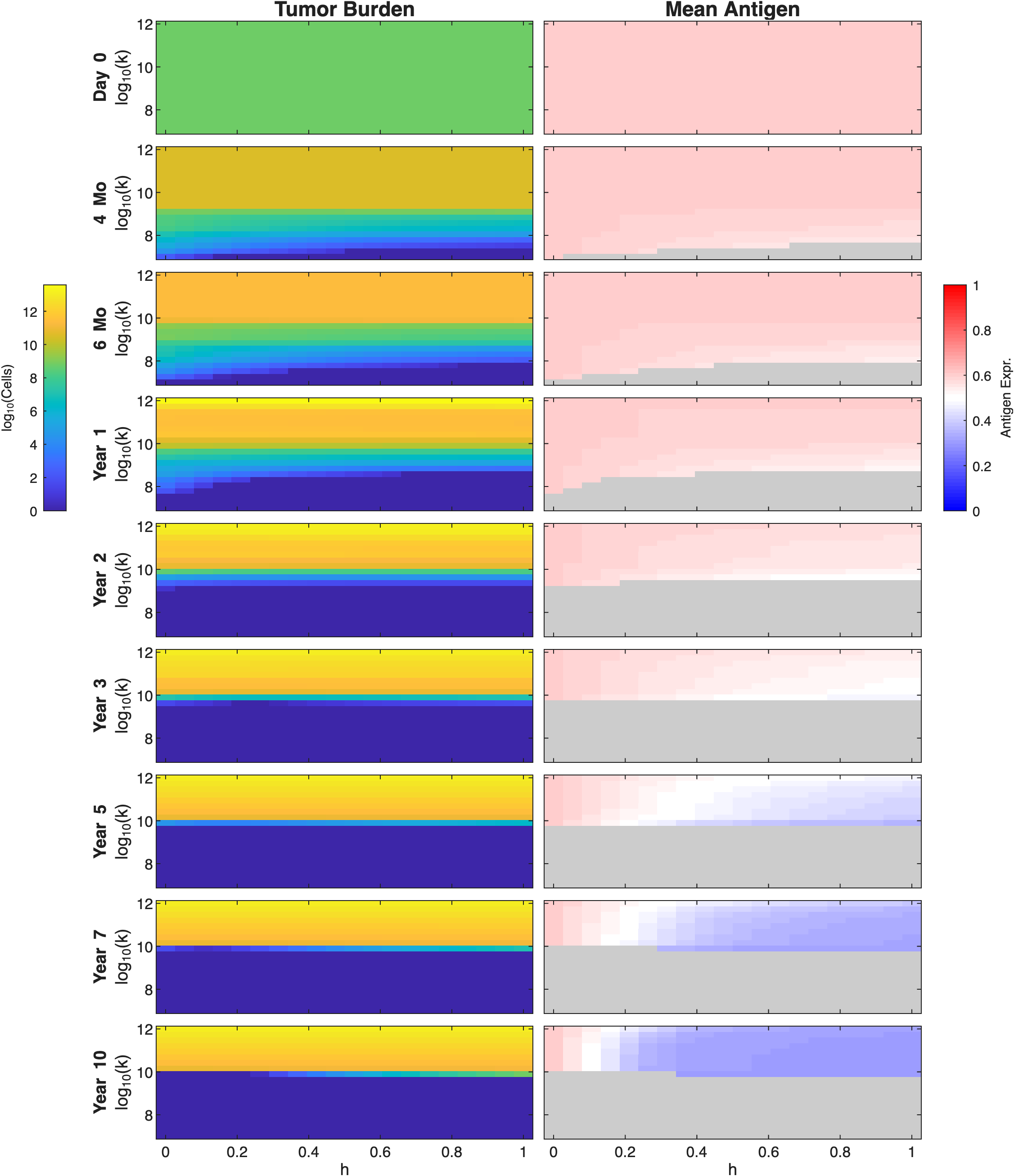}
    \caption{\textbf{Biparametric sensitivity analysis of tumor burden and antigen expression for PDE model from Eqs. \ref{EQ_PDE_model}}. Heatmaps representing the 10-year evolution of a single leukemic clone ($x_0 = 0.6$) under varying CAR T activation thresholds ($k$) and selection coefficients ($h$). The left column displays the total tumor burden ($\log_{10}$ cells), while the right column shows the mean antigen expression (normalized from 0 to 1). Gray areas in the mean antigen plots indicate regions where the tumor burden has been reduced below the detection threshold.}
    \label{fig:antigen_pde_0.6}
\end{figure}

As illustrated in Figure \ref{fig:antigen_pde_0.6}, the system was initialized with a tumor burden of $L_0 = 10^7$ cells using parameters detailed in Table \ref{table param}. We explored a parameter space defined by $k \in [10^7, 10^{12}]$ and $h \in [0, 1]$. Our results indicate a distinct bifurcated behavior dependent on the activation threshold $k$. For a high threshold regime ($k > 10^{10}$), CAR T cells fail to exert sufficient predatory pressure. Consequently, the tumor burden proliferates rapidly, stabilizing at approximately $10^{11}$ cells or higher. In this regime, mean antigen expression remains relatively high ($x > 0.5$) during the initial months, but undergoes an aggressive decline over the 10-year horizon, particularly as $h$ increases.

Conversely, lower values of $k$ ($k < 10^{10}$) enhance a more robust immune response. In these scenarios, higher values of $h$ initially correlate with increased tumor burden during the first few months. However, the system eventually reaches a state where the tumor burden is significantly reduced or stabilized. Notably, for $h$ values approaching 1, the recovery phase is characterized by a significant shift toward lower mean antigen expression, indicative of antigen escape.

These results show that $k$ and $h$ must be considered biparametrically to fully characterize the evolutionary trajectory of the tumor. While a low activation threshold $k$ is essential for controlling tumor volume, the selection parameter $h$ governs the likelihood of immune escape via antigen loss. The optimal therapeutic window for long-term eradication requires a combination of low $k$, thus ensuring high CAR T sensitivity, as well as low $h$ levels to prevent the emergence of aggressive, low-antigen leukemic variants.

\section{Discussion and conclussions}

In this work, we have worked on two models that represent a hierarchical approach to understanding CAR T therapy dynamics. While the Ordinary Differential Equation (ODE) model provides a robust baseline for global population dynamics, the Partial Differential Equation (PDE) model allows us to transcend simple compartments by treating antigen expression as a continuous, evolving landscape. This dual-modeling approach is critical: it showed that while compartmental models are excellent for capturing systemic trends like B cell aplasia \cite{brudno2024car}, they may oversimplify hidden phenotypic shifts that occur within the leukemic population. Our PDE framework captures the subpopulation competition that is often the precursor to clinical relapse, providing a more systematic view of how a tumor escapes the immune system action.

To ensure clinical relevance, we simulated scenarios that mirror the current standard of care and known disease progressions in B ALL, from a controlled CD19$^+$ relapse, to an uncontrolled growth of this type of recurrence or an immune escape to the expression of CD19$^-$ . This included accounting for varying initial tumor burdens, ranging from the high-load environments of initial diagnosis, to the subtle Minimal Residual Disease (MRD) levels \cite{borowitz2015prognostic} between 0.01\% and 10\% of a maximum level of $10^{12}$ cells. By setting initial conditions that reflect pre-lymphodepletion (``bridging therapy''), we simulated the reduction of initial leukemia load $L_0$, showing how this clinical intervention modulates the subsequent CAR T expansion. These simulations demonstrate that the starting patient state, both in terms of tumor volume and the initial distribution of the leukemic clone (in the PDE model, variable $x$), is not just a starting point, but a primary determinant of the therapy's kinetic trajectory over the years. Besides this, we simulated several configurations for the activation threshold $k$, which also describe different qualitative scenarios in terms of possible disease relapses.

Our sensitivity analyses reveal a critical divergence between the two frameworks regarding the primary drivers of tumor control. In the ODE model, long-term outcomes are primarily governed by CAR T persistence $\tau_{CA}$ and the half-maximal activation parameter $h$. However, the PDE model reveals a more stringent bottleneck: the overwhelming dominance of the activation threshold $k$. When antigen expression is heterogeneously distributed, the primary determinant of success is the initial stimulus required to trigger a robust immune response. With $k$ maintaining a relatively high sensitivity index over five years, our results suggest that in patients with low-antigen clones, the therapy's failure is often not due to a lack of ``killing power'', but rather a failure of the CAR T cells to sense the tumor presence at all. The biparametric influence of $h$ and $k$ suggests that patient outcomes are determined early in the treatment timeline. In scenarios where the initial leukemia load $L_0$ is high, the use of "bridging therapy" (pre-lymphodepletion) is essential to modulate the initial tumor burden and prevent immediate CAR T exhaustion or overwhelming Cytokine Release Syndrome (CRS) \cite{morris2022cytokine}. Furthermore, while the phenotypic switching term $\epsilon$ is biologically possible, our sensitivity analysis suggests its impact is secondary to the selection pressure exerted by the $(h, k)$ pair. This reinforces the idea that most relapses are driven by the selection of existing low-antigen clones rather than de novo mutations during the 10-year horizon.

The analysis of leukemic subpopulations highlights a bifurcated escape mechanism. For clones with low antigen expression ($x < 0.5$), the PDE model shows that control is dictated almost exclusively by the activation threshold $k$. These clones effectively achieve immune escape by failing to reach the stimulus threshold necessary for CAR T expansion. Conversely, high-antigen populations ($x \geq 0.5$) are limited by the lifespan $\tau_{CA}$ and expansion rates $\rho_{CA}$ of the effector cells. This suggests that while high-antigen cells are easily targeted, their total eradication is a race against the intrinsic decay of the CAR T population. These findings underscore the importance of considering the MRD levels in clinical settings, as the initial tumor burden significantly dictates whether these thresholds are met. Regarding these last parameters $\tau_{CA}$ and $\rho_{CA}$, a significant finding in our modeling of CAR T dynamics is the trade-off between expansion and safety. In the PDE framework, the expansion rate $\rho_{CA}$ and persistence $\tau_{CA}$ emerge as far more influential than the killing rate $\alpha$. While high $\rho_{CA}$ values ensure long-term surveillance, they also correlate with significantly higher peak populations of active CAR T cells. From a clinical perspective, these high peaks are indicative of a risk for CRS. Our model suggests that therapeutic optimization must find a window where the expansion is sufficient to eliminate low-antigen clones without triggering the hyper-inflammatory responses leading to patient mortality.

The dynamics of healthy B cells also serve as a key agent for CAR T functional persistence. While B cell depletion is initially dominated by the activation threshold $k$, the PDE model shows a temporal shift: by Year 5, the expansion rate $\rho_{CA}$ becomes the primary driver of B cell levels. This indicates that long-term B cell aplasia, a hallmark of successful CAR T therapy \cite{molinos2024impact}, is more closely tied to the sustained expansion of the CAR T population than to the initial activation kinetics. Crucially, our model confirms that the loss of the CAR T population inevitably leads to the reappearance of healthy B cells, providing a clear mathematical link between immune surveillance and clinical off-target effects.

Our analyses presents limitations and possibilities ahead for future work. A inherent limitation of our current deterministic ODE-PDE framework is its behavior during extreme lymphodepletion. When the leukemic population $L(t)$ reaches very low numbers (the MRD tail), stochastic fluctuations (such as the 'extinction probability' of a small clone) become significant. Future iterations of this work could incorporate stochastic differential equations or hybrid agent-based components to better capture the probabilistic nature of total eradication versus the immune escape. Furthermore, while our mathematical framework allows for the inclusion of a baseline population of antigen-negative cells ($L_N(0)$ as an infinitesimal value) we acknowledge that a distinct, pre-existing CD19$^-$ compartment is rarely observed in clinical practice at the time of diagnosis \cite{mikhailova2021relative}. To reconcile this with biological reality, we used a continuous antigen distribution within the PDE model. This approach moves beyond the binary classification of positive or negative expression of the antigen in cells, allowing the emergence of low-antigen phenotypes to arise naturally from the evolutionary dynamics of the system. Future work could benefit from integrating high-resolution single-cell sequencing data to further refine these initial distributions, providing a more granular understanding of whether  the immune escape stems from the expansion of an nearly-invisible population or from a dynamic phenotypic shift across the antigen spectrum. 

Moreover, the present model simplifies CAR T persistence as a linear decay. However, chronic exposure to high tumor burdens often leads to T cell exhaustion, characterized by the upregulation of inhibitory receptors \cite{cillo2024blockade} (e.g., PD-1, LAG-3). Incorporating a state-dependent exhaustion term, where the killing efficiency $\alpha$ and expansion $\rho_{CA}$  diminish as a function of cumulative antigen exposure, would provide a more realistic representation of therapy failure in high-burden patients. Specifically for CAR T cells, a promising future direction is the extension of this PDE framework to dual-targeting CAR T \cite{spiegel2021car} (CD22 and CD19), as modeling the joint distribution of two antigens would allow us to simulate the prevention of antigen escape, which we modulated by the activation threshold $k$. As such, while our PDE model captures phenotypic heterogeneity (antigen expression), it assumes a well-mixed system. In clinical reality, the bone marrow microenvironment provides physical spaces where CAR T infiltration may be limited by stromal interactions or localized immunosuppressive gradients. Transitioning this framework into a spatially-explicit model would allow for the investigation of how the physical architecture of the bone marrow influences the efficacy of CAR T surveillance.

In conclusion, this unified ODE–PDE framework provides several tools for investigating the mechanisms of antigen-dependent relapse in B ALL. Our results emphasize that therapeutic success depends not only on the ``killing efficiency'' of the CAR T cells but, more importantly, on the sensitivity of the immune population to low-density antigen. While this model is currently qualitative and based on literature-derived parameter ranges, it establishes a foundation for future patient-specific calibration. As future work, the transition from a theoretical tool to a clinical decision-support system requires the integration of longitudinal patient data. We could thus focus on using such data to calibrate model parameters in real-time. This \textit{Digital Twin} \cite{niarakis2024immune} approach could allow clinicians to use early CAR T expansion kinetics from the first 14 days post-infusion to predict the 5-year probability of relapse, enabling early intervention with bridging therapies or second-generation infusions. By stratifying patients based on their antigen distribution and activation thresholds, we could then better predict who will achieve long-term remission and who is at risk for antigen-negative escape, leading to better, optimised treatments.

\section{Data availability and supplementary information}
\label{SI}
No datasets were generated or analyzed during the current study.
The source code and functions used in this article can be consulted at \url{https://github.com/salvadorchulian/MathMods-CART-Leuk/}. This repository also includes the full  including different scenarios for different initial conditions and parameters, as well as different leukemic clone configurations. 

The code was executed on a MacBook Pro with the following specifications: 14-inch (2023), equipped with an Apple M3 Pro processor, 18 GB VRAM, running macOS 26.5.1 (25F80).

\section{Acknowdledgements}
The authors thank Manuel Ramírez-Orellana (University Hospital Niño Jesús, Madrid), Álvaro Martínez-Rubio (Marie Curie Institut, Paris) and David Basanta (Moffitt Cancer Center, Tampa, Florida, U.S.A.) for helpful discussions. This work was supported by the Spanish project PID2022-40451OA-I00 funded by Ministerio de Ciencia e Innovación/Agencia Estatal de investigación (doi: 10.13039/501100011033). The support of Asociación Pablo Ugarte (APU, Spain) and Junta de Andalucía (Spain) group FQM-201 is also acknowledged. 

\section{Declaration of generative AI and AI-assisted technologies in the manuscript preparation process}
During the preparation of this work, the authors used Gemini
(powered by Google’s language model, Gemini; \url{https://gemini.google.com/}) 
to improve the readability and language of the work. After using this tool, the authors reviewed and edited the content as needed and take full responsibility for the content of the published article.

\section{CRediT authorship contribution statement}

\textbf{Salvador Chulián}: Conceptualization, Data curation, Formal analysis, Investigation, Methodology, Software, Visualization, Writing – original draft, Writing – review \& editing.
\textbf{Ana Niño-López}: Investigation, Methodology, Validation, Visualization, Writing – review \& editing. 
\textbf{Rocío Picón-González}: Investigation, Visualization, Writing – review \& editing.
\textbf{María Rosa}: Conceptualization, Funding acquisition, Project administration, Supervision,  Writing – review \& editing.
\bibliographystyle{elsarticle-num} 
 \bibliography{bibliography}

\end{document}